\documentstyle[11pt]{article}                     
 \def\Bbb#1{\mbox{\bf $#1$}}                           
\newtheorem{theorem}{Theorem}[section]
\newtheorem{lemma}[theorem]{Lemma}
\newtheorem{proposition}[theorem]{Proposition}
\newtheorem{remark}{Remark}
\newtheorem{remarks}[remark]{Remarks}
\newtheorem{cor}[theorem]{Corollary}
\newcommand{\footer}[1]{{\def\thefootnote{}\footnotetext{#1}}}

\makeatletter
\def\EquationsBySection{\def\theequation{\thesection.\arabic{equation}}
\@addtoreset{equation}{section}}
\makeatother
\EquationsBySection
\renewcommand{\baselinestretch}{1.2}
\def\R{{\Bbb R}}
\def\N{{\Bbb N}}
\def\No{\ov{\N}}
\def\E{\mbox{\sf E}}
\def\P{\mbox{\sf P}}
\def\eps{\varepsilon}
\def\beps{\bm\eps}
\def\<{\langle}
\def\>{\rangle}
\def\Sum{ \displaystyle\sum }
\def\San{\Sum_{\ba\subset[1,n]}}
\def\Frac{ \displaystyle\frac }
\def\Rec#1{ \Frac{1}{#1} }
\def\ex#1{\exp\left\{#1\right\}}
\def\df{\buildrel{\sf df}\over=}
\def\bsT{\mbox{\bf\sf S}}
\def\sA{\mbox{\sf A}}
\def\bsA{\mbox{\bf\sf A}}
\def\bsD{\mbox{\bf\sf D}}
\def\bsV{\mbox{\bf\sf V}}
\def\bsU{\mbox{\bf\sf U}}
\def\sfs{\mbox{\sf s}}
\def\sfc{\mbox{\sf c}}
\def\frac#1#2{{#1\over#2}}
\def\rec#1{ \frac{1}{#1} }
\def\bm#1{\mbox{\boldmath $#1$}}
\def\ca#1{{\cal #1}}
\def\ov#1{\overline{#1}}

\def\refrm#1{{\rm(\ref{#1})}}
\def\set#1{\left\{\,#1\,\right\} }
\def\Proof{{\em Proof. \hspace{10pt} }}
\def\otk{{\otimes k}}

\def\bx{{\bm X}}
\def\BX{{\Bbb X}}
\def\L{{\Bbb L}}
\def\BE{{\Bbb E}}

\def\l({\bm{\big(\hspace{-3pt} \big(}\, }
\def\r){\,\bm{\big)\hspace{-3pt} \big)} }
\def\ll{\Big\langle\, }
\def\rr{\,\Big\rangle}
\def\lll{\,\bm{\Big\langle \hspace{-4.5pt} \Big\langle}\, }
\def\rrr{\,\bm{\Big\rangle \hspace{-4.5pt} \Big\rangle}\, }

\def\wh#1{\widehat{#1}}

\def\QED{\hfill \rule{2.5mm}{2.5mm}
\vspace{7pt}

}
\def\Int{\displaystyle\int }
\def\Ib#1{1\hspace{-3pt} {\rm I}_{\set{#1}}}
\def\wt#1{\widetilde{#1}}
\def\iid{independent identically distributed }
\def\ba{{\bm\alpha}}
\def\bb{{\bm\beta}}
\def\fba{f_{\ba}}

\def\imp{~ \Rightarrow~ }

\def\qand{\quad\mbox{and}\quad}
\def\ttN{\set{0,1}^\N}
\def\LNM{Lecture Notes in Math.}
\def\CWRU{Case Western Reserve University}
\def\PMS{Probab.\ Math.\ Statist.}
\def\PTRF{Probab.\ Theory Related Fields}
\def\AnMS{Annals of Math.\ Statist.}
\def\AP{Ann. Probab.}
\def\AU{Auburn University}
\def\AJM{Amer.\ J.\ Math.}
\def\PAMS{Proc.\ of Amer.\ Math.\ Soc.}
\def\SM{Studia Math.}
\begin{document}
\footer{{\bf AMS (1980) Classification} {\em  Primary}: 60B11, 60H07, 15A69
{\em Secondary}: 46M05 , 46E30, 60E15, 62H05,  62G30, 42A55. 15A52, 26D15}
\vspace{10pt}

{\Large \bf Infinite order decoupling 
of random chaoses in Banach space}
\vspace{5pt}

{{\bf Jerzy Szulga}} \footnote{Dept. of Mathematics, Auburn University, Auburn,
AL 36849}
\vspace{7pt}

\hfill  \rule{5.4in}{1pt}
\vspace{5pt}

\hfill  \parbox{5.4in}{\small We prove a number of decoupling inequalities for
nonhomogeneous random polynomials with coefficients in  Banach space. Degrees
of homogeneous components enter into comparison as exponents of multipliers of
terms of certain Poincar\'e-type polynomials. It turns out that the fulfillment
of most of types of decoupling inequalities may depend on the geometry of
Banach space.}
\vspace{5pt}

\hfill  \rule{5.4in}{1pt}
\vspace{5pt}

\hfill  \parbox{5.4in}{{\bf \small KEY WORDS}: \small decoupling principle,
symmetric tensor products, random polynomials, multiple random series, multiple
stochastic integrals, random multilinear forms, random chaos, Gaussian chaos,
Rademacher chaos, stable chaos, multiple Wiener integral, multiple stable
integral, Mazur-Orlicz polarization formula, symmetrization, Banach space,
Banach lattice, Krivine's type, rearrengement invariant space, convexity,
Orlicz space, Rademacher sequence, Gaussian law, Walsh polynomials, empirical
measure.}
\vspace{3pt}

\begin{center}
\rule{6in}{1pt}
\end{center}
\vspace{7pt}

\renewcommand{\baselinestretch}{0.8}
{\small\sf \tableofcontents}
\renewcommand{\baselinestretch}{1.1}
\eject
\section{INTRODUCTION}
The concept of decoupling stems from the martingale theory (cf. the survey
\cite{Bur:mf}).
First decoupling inequalities for multiple random forms were proved in
\cite{McCT:realdec,McCT:bandec,Kwa:dec}, and many variants have been published
since the time when the above papers were published (to name but a sample, cf.
\cite
{DeA:dec,Hit:tan,Zin:comp,NolP,Pen,KalS,Szu:Umli,RosW:clock,RST,Szu:dechom},
and further references in there).  So far, all known results have involved a
two-sided estimate of $L^p$-(or Orlicz) norms of suitable $k$-homogeneous
multilinear forms (or multiple integrals), where $k$ is arbitrary but fixed.
Decoupling constants are degree-dependent and escape to infinity. If the degree
increases, the strength of the decoupling principle seems to decline.

In this paper we will show how to overcome this deficiency  (one cannot expect
that decoupling constants remain bounded). Our approach is based on a suitable
normalization of polynomials
\[
Q(\BX;t)=Q_0+t Q_1(\BX) +\ldots +t^n Q_n(\BX),
\]
where $\BX=[X_{ij}]$ is a matrix of random variables, with independent rows,
and $Q_k$ is a Banach space-valued homogeneous polynomial of degree $k$ (a
$k$-linear form). Under suitable integrability and
symmetry assumptions the presented
decoupling principle compares norms (e.g.,
Orlicz norms) of $Q(\BX,t)$ and of the polynomial $Q(\ov{\BX})$,
\[
\|Q(\BX)\|\sim\|Q(\ov{\BX})\|
\]
where the matrix $\ov{\BX}$ is a ``decoupled'' version of $\BX$, i.e., the
columns of $\BX$ are replaced by their independent copies. For example,  on the
real line, for Rademacher or standard Gaussian random variables, we check
directly that
\begin{equation}\label{realdec}
\E|\sum_{k\ge 0}  Q_k|^2 =
\E|\sum_{k\ge 0}\overline{Q}_k|^2
\end{equation}
where
\[
\overline{Q_k}=\left\{
\begin{array}{ll}
\Rec{\sqrt{k!}}\Sum_{i_1,\ldots,i_k, i_j\neq i_{j'}} f_k(i_1,\ldots,i_k)
X_{1i_1}\cdots X_{ki_k},& \mbox{if $f_k$ is symmetric, or}\\
\Sum_{i_1<\ldots<i_k} f_k(i_1,\ldots,i_k) X_{1i_1}\cdots X_{ki_k},&\mbox{if
$f_k$ is tetrahedral}\\
\end{array}
\right.
\]
and $Q_k$ follows the same pattern, respectively, but without the multiplier
$1/\sqrt{k!}$ in the symmetric case. Also, $L^2$-norm can be replaced by
$L^p$-norm, $p>1$, at cost of multiplying each $k$-homogeneous polynomial by a
constant $c_p^k$.

Degrees of specific components will enter into formulas as exponents of certain
multipliers. We replace, so to speak, external constants by internal constants
or, more precisely, by sequences of constants. Asymptotic behavior of such
sequences is of interest, and the exponential growth is most desirable.  If
columns of $\BX$ are identical, we call $Q(\BX)$ a `` random chaos'', and when
they are independent (desirably -- identically distributed), a ``decoupled
random chaos''. A tetrahedral decoupled chaos can be written as a sort of
lacunary chaos, by a monotone (non-unique, in general) change of ordering on
all tetrahedra, from the coordinatewise ordering to a linear ordering. An
analogous procedure for symmetric chaoses is possible ``locally'', i.e., for a
fixed and finite  order, and when only a finite number of random variables is
involved.

Thus, decoupling inequalities can be viewed as embedding-projection theorems.
In the infinite order decoupling we require that projections are contractions.
We will observe new phenomena, absent in the homogeneous (or finite order)
case.  First of all, in the infinite order decoupling, two types of
inequalities (``the lower decoupling'' -- domination by the chaos, and ``the
upper decoupling'' -- domination of the chaos) determine two distinct problems.
Already homogeneous tetrahedral and symmetric chaoses behave differently (cf.
Bourgain's example in \cite{McCT:bandec}).

A robust lower decoupling inequality is satisfied, i.e., the inequality is
fulfilled in any Banach space and for an arbitrary symmetric integrable
polynomial chaos. At the same time, the fulfillment of a robust upper
decoupling inequality is still uncertain.  In order to study the upper
decoupling principle, we introduce a class of Banach spaces that are
characterized by a certain inequality involving linear forms in independent
random variables with vector coefficients (in some aspects, the property is
similar to classical properties of Banach spaces, like Rademacher type and
cotype, smoothness or convexity of norm, etc.

The main feature of the new class of spaces is that random polynomials admit a
``horizontal slicing'', reducing the study to that of sums of independent
random variables.  In the introduced class of Banach spaces a sign-randomized
upper decoupling inequality holds.  The class of spaces allowing slicing of
random polynomials  is, unfortunately, geometry-dependent, and very fragile.
It is sensitive to an equivalent renorming (just an addition of a
two-dimensional normed space with $\ell^\infty$- or $\ell^1$- norm terminates
the property), in contrast to the homogeneous case. Therefore, a positive
result will always require the existence of a suitable equivalent norm (cf.
``smooth'' vs. ``smoothable'', or ``convex'' vs.  ``convexifiable'').  The
family of Banach spaces, satisfying the slicing requirements, contains Banach
lattices of finite cotype and sufficiently convex norm.  This class turns out
to be suitable even for the tetrahedral lower decoupling.

In Section 2 we introduce the nonhomogeneous tensor product notation, define
the domination, and derive some basic relations. We refer to \cite{Per} and the
literature included there for a treatment of Gaussian symmetric tensors. The
comparison will be given in terms of the aforementioned Poincar\'e -type
polynomials or, equivalently, in terms of a semigroup of contractions,
associated with a random polynomial.

The new results are gathered in Sections 3 and 4.  The employed techniques in
the nonhomogeneous case are different from techniques related to the
homogeneous case. First of all, the applicability of conditional expectations
is limited. Secondly, the type of domination forbids use of external constants.
We will use a ``slicing technique'', which reduces the study to a case of
certain sums of random variables.

In Section 5 we will indicate some directions in a further study of decoupling
inequalities, and show that many assertions can be directly obtained from
results of this paper. For example, one can formulate decoupling inequalities
in the language of infinite order stochastic multiple integrals.  We will also
collect some observations that do not fit into the main line of the paper,
although may be of some interest.

Let us point out that a widely understood convexity is the principal feature
implicit in most of applied techniques. This includes the setting of Banach
spaces and existence of moments (integrability) of involved random variables.
One can find a number of decoupling principles in the literature, where the
convexity is of no concern, and the focus is on positivity, not on symmetry.
However, most of the known results have been obtained so far at the cost of the
limitation to the real line (see \cite{KalS} for a discussion on the latter
topic). In some special cases, e.g., for Gaussian homogeneous polynomials, a
decoupling principle applies to probabilities $\P(\cdot \notin K)$, where $K$
is a convex symmetric set in a Banach space \cite{Kwa:dec}; the aforementioned
paper \cite{DeA:dec} deals with $p$-stable random variables and spaces $L^r$,
$0<r<p$, etc.  

\section{RANDOM TENSOR PRODUCTS}
\def\BY{{\Bbb Y}}

\subsection{Notation}
Throughout the paper $\beps=(\eps_i)$ denotes a Rademacher sequence, that is,
$\eps_i$ are independent random variables taking values $\pm 1$ with
probability $1/2$. Walsh functions are products of Rademacher random variables.
We will denote by $\bx$ a sequence, and by $\BX=[\bx_1,\ldots,\bx_n]$  a
matrix, of real random variables.  $(\BE,\|\cdot\|_\BE)$ denotes a real Banach
space. By $(\L,\|\cdot\|_\L)$  we denote a rearrangement invariant Banach space
of integrable random variables, $\L\subset L^1(\P)$, defined on a probability
space $(\Omega, \ca F, \P)$, rich enough to carry independent sequences, and
with a separable $\sigma$-field. In fact, we will use only specific properties
of $\L$, ensured by the above restrictions.
\[
\begin{array}{ll}
\mbox{\hspace{-10pt} \bf (L)}\hspace{20pt} &
\mbox{\parbox{5in}{\em 
Conditional expectations are contractions acting in $\L$}}\\
\end{array}
\]
$\L(\BE)$ denotes the Banach space of $\BE$-valued random variables (i.e.,
strongly measurable mappings from $\Omega$ into $\BE$) whose norms belong to
$\L$ a.s., and let $\|\theta\|_{\L(\BE)}=\|\,\|\theta\|_{\BE}\,\|_{\L}$.
Whenever it causes no ambiguity, we omit the subscript.

Let $\N=\set{1,2,\ldots}$ be the set of natural numbers, and $\No=\N\cup
\set{0}$. For $m,n\in\No$, put $[m,n]=\set{m,m+1,\ldots,n}$. 
Throughout the paper, the bold-face Greek characters $\ba, \bb, \ldots $, etc.,
will denote subsets of $\N$, identified with $\set{0,1}$-valued sequences:
\[
\N\supset \ba ~ \longleftrightarrow~
\ba=(\alpha_1,\alpha_2,\ldots)\in\set{0,1}^\N.
\]
Denote $|\ba|=\# \ba=\sum_i\alpha_i$, and $\ba'=(1-\alpha_1,1-\alpha_2,
\ldots)$.  The following convention will be very helpful.  Let $\BE$ be a
nonvoid set.  Suppose that $0,1\in \BE$. Define two operations $\set{0,1}\times
\BE\to \BE$:
\[ 
0x=0,\quad 1x=x\qand x^0=1,\quad x^1=x \mbox{ (by convention, }0^0=1).
\] 
 If $*$ is any operation in an abstract set $Z$, then we will use the same for
functions taking values in $Z$. In particular, if $*:\BE\times Y\to Z$, then,
by writing $*:\BE^\N\times Y^\N\to Z^\N$, we understand the action of $*$
coordinatewise. For example, $(\bx * \bm y)_i=x_i * y_i$, $i\in \N,
\bx=(x_i),\,\bm y=(y_i)$. Also, for $\bx\in \BE^\N$ and $\ba\in
\ttN$, $\ba \bx=(\alpha_i x_i)$ and $\bx^\ba=(x_i^{\alpha_i})$. 
In section 3.4, the term $\Bbb S X$, where $\Bbb S$ and $\BX$ are $(N\times n)$
matrices, according to our convention, will denote a new $(N\times n)$ matrix,
whose entries are products of entries of $\Bbb S$ and $\Bbb X$.

Identifying $\alpha$ and  a constant sequence $(\alpha,\alpha,\ldots)$, we have
then $\alpha \bx=(\alpha x_1,\alpha x_2,\ldots)$.  For $\ba \subset \N$, we
identify $\N^\ba$ and the subset $\set{\ba \bm i:\bm i \in
\N}$ of $\No^\N$ (the empty set is identified with $\set{(0,0,0,\ldots)}$. 
We will consider functions $\bm f=(f_\ba:\ba\in \ttN)$, where $f_\ba:\N^\ba\to
\BE$ ($f_\emptyset$ is an element of $\BE$, and, if necessary, $\bm f $ may be
identified with suitable functions $\bm f:2^\N\times \No^\N\to \BE, \,\bm
f(\ba,\bm i)=f_\ba(\ba \bm i)$). For definitness, we request that all functions
$\bm f =(\fba)$ under consideration have a finite support (i.e. $f(\ba,\bm
i)=0$ for all but finitely many $\ba$ and $\bm i$.

In this paper we will see an abundance of summation, averaging, and integration
on several levels, in order to diminish the notational burden, we introduce a
variety of summing brackets  (of course, we might replace all following
brackets by just one but, by doing so, we would cause a serious visual
dissonance).  Define
\[
\ll \fba \rr =\ll \fba \rr _\ba=\sum_{\ba\bm i}\fba(\ba\bm i)\qand\lll f \rrr
=\sum_\ba \ll \fba \rr _\ba.
\]
and, for $\bm x = (x_i) \in \BE^\N$,
\[
\l( \bm x \r) =\sum_i x_i
\]
All functions $\bm f=(f_\ba)$, appearing in the sequel are assumed to vanish on
diagonals, i.e., $f_\ba(\ba \bm i)=0$ if at least two nonzero-arguments $(\ba
\bm i)_k$ are equal.
Define the symmetrizator $\wh{\bm f}$, which unifies values of functions $\fba$
on distinct tetrahedra, by the formula
\[
\wh{\fba}=\rec{|\ba|!}\sum_\sigma \sigma \fba,
\]
where the sum is taken over all permutations $\sigma$ of $\ba$, and
$\sigma\fba=\fba\circ\sigma$. If $\wh{\bm f}=\bm f$, then the functions is
called symmetric.  Call a function $\bm f$ {\em tetrahedral}, if it may take
nonzero values only on the main tetrahedron:  $\ba=\set{1,\ldots,|\ba|}$,
$i_1<\ldots<i_{|\ba|}$ .

The random matrix $\BX=(X_{ij}:i\in [1,n],j\in \N)
\in \R^{[1,n]\times \N}$, considered before, can be written as
$\BX=[\bx_1,\bx_2,\ldots]\in (\R^\N)^{[1,n]}$, where $\bx_i=(X_{ij}:j\in\N)$.
Define a {\em tensor product } $\BX^\otimes=(\BX^{\otimes\ba})$ on
$\R^{[1,n]\times\N}$ by the formula
\[
\BX^{\otimes\ba}(\ba\bm i)=X^{\alpha_1}_{1\,i_1}\cdots X^{\alpha_n}_{n\,i_n}, 
\]
and the {\em symmetric tensor product}, by
$\bx^{\wh{\otimes}}=\wh{\bx^\otimes}$.  By convention, a single sequence $\bx$
can be viewed as a matrix (the sequence) $[\bx,\ldots,\bx]$. Whence the tensor
products $\bx^\otimes$ and its symmetrizations are well defined. One can
consider other type of symmetrizators $\sf U$  and the induced symmetric tensor
products ${\sf U}\otimes$ (see Section 5).
\vspace{5pt}

The Mazur-Orlicz polarization formula can be written as follows (it is
fulfilled in  any commutative algebra).

\begin{equation}\label{MOalpha}
\begin{array}{rl}
\bx^{\wh{\otimes}\ba}
=&\displaystyle\rec{k!}\Sum_{\bb\subset \ba}
(-1)^{k-|\bb|}\l(\bb\bx\r)^{\otimes\ba},\\ =&\displaystyle
2^{-n}\rec{k!}\Sum_{\bb\subset [1,n]}
(-1)^{k-|\bb\ba|}\l(2\bb\ba\bx\r)^{\otimes\ba},\\
\end{array}
\end{equation}
where $|\ba|=k$.  Endowing $2^{[1,n]}$ with the uniform probability, the
functions
\[
r_i(\bb)=(-1)^{\beta_i},\qquad \bb\in 2^{[1,n]},\quad i=1,\ldots,n
\]
are representations of the first $n$ Rademacher random variables.  We define
{\em Walsh functions} as products of Rademacher functions:
\begin{equation}\label{walsh}
w_\ba(\bb)\df \prod_{i\in \ba} r_i(\bb)=(-1)^{\otimes \ba\bb},\qquad\bb\in
2^{[1,n]}.
\end{equation}
Notice the presence of Walsh functions in the Mazur-Orlicz polarization
formula.

By $\ca S^0_k=\ca S^0_k(\BX)$ denote the $\sigma$-field generated by the family
\[
\set{g_\ba(\BX^\ba)\,:\, \ba\subset [1,n],\,|\ba|=k,\,
g_\ba=\hat{g_\ba},\, g:\R^k\to \R},
\]
and let $\ca S_k=\ca S_k(\BX)$ be  the $\sigma$-field generated by $\ca
S_1^0\cup\ldots\cup \ca S^0_k$, and $\ca S$ be spanned by $\bigcup_k \ca S_k$.

\begin{proposition}\label{condextensor}\rm
Let $[\bx,\BX]$ have independent rows and i.i.d columns. Let $\bm g_\ba$ be a
function vanishing on diagonals and equal 1 off diagonals. Then the following
equalities hold.  %
\begin{equation}\label{bax}
g_\ba \l(\bb\ba\BX\r)^{\otimes\ba}=
\E[\,g_\ba \l(\ba\BX\r)^{\otimes\ba} \,|\,\BX^\bb \,];
\end{equation}
\begin{equation}\label{BXbx}
g_\ba (\BX+\BX')^{\otimes\ba}=
\E[\,g_\ba (2\BX)^{\otimes\ba}\,|\,\BX+\BX'\,],
\end{equation}
where $\BX'$ is an independent copy of $\BX$; %
\begin{equation}\label{BXbx1}
g_\ba \BX^{\otimes\ba}=\E[g_\ba \,(\BX+\BX')^{\otimes\ba}\,|\,\BX\,]
\end{equation}
provided $\E\BX'=0$, and $\BX$ and $\BX'$ are independent; %

\begin{equation}\label{Vcontr}
\lll \bm f \BX^{\wh{\otimes}} \rrr =
\lll \bm f \E[\BX^{\otimes}|\ca S(\BX)] \rrr =
\E[\lll \bm f \BX^{\otimes}\rrr|\ca S(\BX) ] 
\end{equation}
\end{proposition}                                         
\Proof
Conditions \refrm{bax}, \refrm{BXbx}, and \refrm{BXbx1}  follow immediately.
In order to see the fulfillment of the remaining condition, it suffices to
implement the following simple rule.  For two $\sigma$-fields $\ca F_1,\,\ca
F_2$, and  random variables $Z_m$, $m=1,\ldots,M$, if $\E[\,Z_m\,|\,\ca
F_1\,]=\E[\,Z_1\,|\,\ca F_1\,]$, and $\sum_mZ_m$ is $\ca F_1$-measurable, then
$\E[\,Z_i\,|\,\sigma(\ca F_1\cup\ca F_2]=\sum_m Z_i/M$.  This rule allows us to
reduce each situation to the homogeneous case, and then the proof is
direct.\QED

\subsection{Domination of random polynomials}
The decoupling principle for nonhomogeneous polynomials will be defined in
terms of a more general concept of domination, applied to certain
Poincar\'e-type polynomials (e.g., cf. various variants of hypercontraction
\cite{Gro,KraS:hyp,KwaS} or Malliavin's calculus, cf. \cite{Sug}).

First, we decide a setting of domination.  Let $(\BE,\|\cdot\|)$ be a Banach
space, $\L\subset L^0(\P)$ be an algebra of re\-al ran\-dom va\-ria\-bles,
en\-dowed with a po\-si\-tive functional $\varphi$, and $\L(\BE)=\set{\xi\in
L^0(\BE):\|\xi\|\in\L}$ a.s. Define the functional
$\Phi(\xi)=\varphi(~\|\xi\|~)$.

Usual examples consist of $L^p$-norms or quasi-norms, $0\le p\le\infty$, Orlicz
(or more general rearrangement invariant) norms, distribution tails
$\phi(\theta;t)=\P(|\theta|>t)$, etc. (cf. \cite[Chapters 3, 5]{KwaW:book} for
more examples).

Consider   $\BE$-valued random polynomials $\lll \bm f \BX^\otimes \rrr $,
where $\bm f =(\fba:\ba\subset\N$, $\quad \fba:\N^\ba\to \BE$, $f_\ba\equiv 0$
for all but finitely many finite sets $\ba$, and $\bm f$ belong to a certain
category $\ca F$ of functions, realized on the class of Banach spaces.  The
role of a constant is to be played by a real valued function $\bm c
=(c_\ba):2^{\N}\times \No\to\R\in\ca U$, where $\ca U$ is realized on $\R$.
The system $\ca S=(\BE, \L, \varphi, \ca F)$ becomes the  setting of
domination.

{\bf Definition}. ~ Say that $\BX$ is {\em dominated} by $\BX'$ in the setting
$\ca S $ with a constant $\bm c$ ($\BX\preceq_{\ca S,\bm c} \BX'$, in short),
if, for every $n\ge 0$,
\[
\Phi(\lll \bm f \BX^\otimes \rrr )\le
\Phi(\lll \bm f\bm c \BX'^\otimes \rrr),\qquad \bm f \in\ca F.
\]
If the constant is of the form $c_\ba =c^{|\ba|}$, where $c$ is a positive
number we will say that the $\BX$ is {\em exponentially dominated} by $\BX'$.

It is easy to see that linear forms in zero mean integrable random variables
are comparable with their symmetrized counterparts. A similar property is
enjoyed by random chaoses. Let $\Phi$ be a positive convex functional defined
on $L^1(\BE)$ turning conditional expectations into contractions, i.e.,
\begin{equation}\label{Phi}
\Phi(\E[ X|\ca F])\le \Phi(X)
\end{equation}
(for example, $\Phi(X)=\|X\|_{\L(\BE)}$, where $\L$ is a rearrangement
invariant space of real random variables, or $\Phi(X)=\E\varphi(\|X\|_\BE)$,
where $\varphi$ is an increasing convex function).
\begin{proposition}\label{dec:l:sym}\rm
Let   $\BX$ and $\BX'$ be independent identically distributed matrices  with
independent rows and interchangeable columns (this includes the case of
matrices with identical columns). let $\bm f =(f_\ba)$ be an $\BE$-valued
function vanishing on diagonals.
\begin{enumerate}
\item 
Then
\begin{equation}\label{dec:e:sym}
\Phi(\lll \bm f(\BX-\E\BX)^{\otimes} \rrr)\le 
\Phi(\lll \bm f(\BX-\BX')^{\otimes} \rrr).
\end{equation}
\item 
If $\bm f$ is  symmetric, and $\BX$ has independent columns, then there exists
a Walsh system $w_\ba$, independent of $\BX$ and $\BX'$, such that
\begin{equation}\label{dec:e:zeromean}
\Phi(\lll \bm f(\BX-\BX')^{\otimes} \rrr )\le 
\Phi(\lll \bm f\bm w (2\BX)^{\otimes} \rrr ). 
\end{equation}
where $(\bm f\bm w)_\ba(\bm i)=\fba(\bm i) w_{\ba}$.
\item
If $\bm f$ is symmetric, and $\E\bx=0$, then
\begin{equation}\label{coupdec:e:zeromean}
\Phi(\lll \bm f(\bx-\bx')^{\otimes} \rrr )\le 
\Phi(\lll \bm f\bm w (4\bx)^{\otimes} \rrr ). 
\end{equation}
\end{enumerate}
\end{proposition}
\Proof 
Inequality \refrm{dec:e:sym} follows by convexity and contractivity of
conditional expectations.

Inequality  \refrm{dec:e:zeromean}  follows from the estimates
\[
\begin{array}{rl}
\Phi(\lll \bm f{\BX-\BX'}^{\otimes} \rrr)
&=\Phi(\San
\ll \bm \fba (\BX-\BX')^{\otimes\ba} \rr )\\
&=\Phi(\San
\ll \bm \fba \Sum_{\bb\le\ba}\BX^{\otimes\bb}\otimes
(-\BX')^{\otimes(\ba-\bb)} \rr )\\
&=\Phi(\San
\ll \bm \fba 
2^{|\ba|-n}\Sum_{\bb\subset[1,n]}\BX^{\otimes\bb\ba}\otimes
(-\BX')^{\otimes(\ba\bb')} \rr )\\
&\le\rec{2^n}\Sum_{\bb\subset[1,n]}
\Phi(\San \ll \bm \fba 
2^{|\ba|}\BX^{\otimes\bb\ba}\otimes
(\BX')^{\otimes(\ba\bb')}(-1)^{\otimes(\ba\bb')} \rr )\\
&=\rec{2^n}\Sum_{\bb\subset[1,n]}
\Phi(\San w_\ba(\bb)\ll \bm \fba
(2\BX)^{\otimes\ba} \rr )\\
&=\Phi(\San\ll \bm \fba w_\ba
(2\BX)^{\otimes\ba} \rr ).\\
\end{array}
\]
The proof of \refrm{coupdec:e:zeromean} is similar. However, the assumption
$E\bx=0$ is essential in the following argument. We have

\begin{eqnarray*}
\Phi(\lll \bm f(\bx-\bx')^{\otimes} \rrr )
 =\Phi(\San
\ll \bm \fba (\bx-\bx')^{\otimes\ba} \rr )\\
 =\Phi(\San
\ll \bm \fba \Sum_{\bb\le\ba}\bx^{\otimes\bb}\otimes
(-\bx')^{\otimes(\ba-\bb)} \rr )\\
 =\Phi(\San
\ll \bm \fba 
2^{|\ba|-n}\Sum_{\bb\subset[1,n]}\bx^{\otimes\bb\ba}\otimes
(-\bx')^{\otimes(\ba\bb')} \rr )\\
 \le\rec{2^n}\Sum_{\bb\subset[1,n]}
\Phi(\San \ll \bm \fba 
2^{|\ba|}\bx^{\otimes\bb\ba}\otimes
(\bx')^{\otimes(\ba\bb')}(-1)^{\otimes(\ba\bb')} \rr )\\
 \le\rec{2^n}\Sum_{\bb\subset[1,n]}
\Phi(\San \ll \bm \fba 
2^{|\ba|}(\bx+\bx')^{\otimes\bb\ba}\otimes
(\bx+\bx')^{\otimes(\ba\bb')}(-1)^{\otimes(\ba\bb')} \rr )\\
 =\rec{2^n}\Sum_{\bb\subset[1,n]}
\Phi(\San w_\ba(\bb)2^{|\ba|}\ll \bm \fba 
\E[(2\bx)^{\otimes\ba}|\bx+\bx'] \rr )\\
 \le\Phi(\San\ll \bm \fba w_\ba
(4\bx)^{\otimes\ba} \rr).\\
\end{eqnarray*}
The proof has been completed. \QED

The following {\em contraction principle} is well known in the one dimensional
case (cf. \cite{Kah} for the real case, and \cite{Hof}, for the vector case).
\begin{theorem}\label{multicontr}\rm
Let $\varphi:\R_+\to\R_+$ be a convex increasing function. Let $\BX$ be a
matrix of real symmetric random variables with independent rows and either
independent or identical columns. Then for every $\BE$-valued function $\bm
f\in \ca F_S$ (or $\ca F_T$ ), and bounded real function $\bm g=(g_\ba)$ with
$c=\|g\|_\infty$, of the form $g_k(\bm i)=g_{k1}(i_1)\cdots g_{kk}(i_k)$, we
have 
\[
\E\varphi\left(\|\lll \bm f\bm g \BX^\otimes \rrr\|\right)\le
\E\varphi\left(\|\lll \bm f  (c\BX)^\otimes \rrr\|\right)
\]
\end{theorem}
\Proof
In the case of a homogeneous  chaos, i.e., when columns of the matrix $\BX$ are
identical, the result appeared in \cite{Kwa:dec}, while for nonhomogeneous 
chaoses, in \cite{KraS:hyp}.
The case with independent columns follows by a spreading argument. That is,
there exists an enumeration of functions $\bm f$ such that the polynomial $\lll
\bm f \BX^\otimes\rrr$ can be written as a polynomial $\lll \bm f'
\bx^\otimes\rrr$. Hence, we arrive in the previous situation. \QED

Say that tails of two random variable $X$ and $Y$ are comparable, if, for some
constants $K>0$ and $t_0\ge 0$
\[
\P(|X|>t)\le K\P(|Y|> K t)\qand
\P(|Y|>t)\le K\P(|X|> K t), \qquad t\ge t_0.
\]
Notice that, at cost of increasing the constant $K$, one may assume that the
above estimates hold for every $t>0$. Thus, there exist  probability spaces and
copies $X',X''$ and $Y',Y''$ of $X$ and $Y$, respectively, such that $|X'|\le
K' |Y'|$ and $|Y''|\le K'|X''|$ a.s. In particular, the upper decoupling
inequality is satisfied simultaneously for chaoses spanned by $\bx$ and $\bm
Y$, provided components of both sequences are independent and have comparable
tails.

\begin{cor}\label{comp}\rm
Let $\BX_1$ and $\BX_2$ be matrices of real symmetric random variables with
independent rows and independent or identical columns. Suppose that
corresponding entries of both matrices have comparable tails. Then, for any
symmetric or tetrahedral function $\bm f$, polynomials in $\BX_1$ and $\BX_2$
are comparable, i.e., for any increasing function $\phi:\R_+\to \R_+$,
\[
\E\varphi\left(\|\lll \bm f \BX_i^\otimes \rrr\|\right)\le
\E\varphi\left(\|\lll \bm f (c\BX_j)^\otimes \rrr\|\right),
\]
$i,j\in\{1,2\}$, for some constant $c$, depending on the tail domination
constant $K$.
\end{cor}

\subsection{Lower and upper decoupling inequalities}
{\bf Definition.} Let $\bx$ be a sequence of real independent random variables
and $\BX$ be a matrix whose columns are independent copies of $\bx$.  Let $\ca
F$ be a class of functions $\fba=(f_\ba)$, $f_\ba:\N^\ba\to\BE$. Denote by
$\ca{UD}=\ca{ UD} (\BE;\Phi;\ca F)$ (respectively, $\ca{LD}=\ca{LD}
(\BE;\Phi;\ca F)$) the class of sequences $\bm X=(X_i)$ of independent random
variables (more exactly, the class of product probability measures) such that
that the {\em upper decoupling inequality} (respectively, the {\em lower
decoupling inequality}) holds on $\ca F$, i.e., there exists a constant $c$
such that, for every $n\in\N$ and $\bm f\in \ca F$, one has
\[
\E\Phi(\lll \bm f \bx^\otimes \rrr)\le
\E\Phi(\lll \bm f (c \BX^\otimes \rrr)
\qquad\mbox{(respectively,}\quad  
\E\Phi(\lll \bm f (\BX^\otimes \rrr)\le
\E\Phi(\lll \bm f (c\bx)^\otimes \rrr)\mbox{ ).}
\]
If the considered sequences have components with the same probability
distribution $\mu$, we will say that $\mu$ (or a random variable with the
distribution $\mu$) satisfies the upper (respectively, lower) decoupling
inequality.

The most important classes are ${\ca F}_S$, the class of  symmetric functions,
and ${\ca F}_T$, the class of tetrahedral functions (recall that we always
assume that functions vanish on diagonals). One can consider also other classes
(cf. Section 5).
More precisely, the decoupling introduced above is understood in the sense of
the exponential domination. Note that in most cases of interest that is a
desired property. By the same token one can discuss the decoupling in a weaker
sense (with a ``constant'' $c_\ba$ being not necessarily of the exponential
type), but then both sides of decoupling, the lower and upper inequality,
should be treated separately.

Proposition \ref{dec:l:sym} indicated that in case of insufficient symmetry it
is necessary to randomize signs of consecutive homogeneous components of a
random chaos. In the proposition, such a randomization does not affect internal
components of homogeneous polynomials. However, as will be shown, frequently
one needs random signs within each and every homogeneous term, and the
intricacy of such a randomization may vary. for example,  one may use Walsh
multipliers induced either by one sequence $\beps$, or by a matrix
$[\beps_1,\beps_2,\ldots]$ with independent Rademacher columns, or, instead of
Walsh functions, one may require plain Rademacher family, indexed by the
multi-index $\ba \bm i$. The latter sign-randomization is the only known way,
so far, of extending Theorem \ref{multicontr} to functional multipliers $\bm g$
whose arguments are not separated (cf. \cite{KraS:sum}).


\section{SLICING AND DECOUPLING}
\subsection{Slicing}
Let $\ca C$ be a class of finite random real sequences and $N$ be an integer.
Say that an $n\times N$ random real matrix  $\BX$ is {\em $\ca C$-sliceable},
if its rows are independent and belong to $\ca C$.
For $\ba \subset [1,N]$, denote
$\ba^*=\max\ba =\max\set{i\in[1,N]:\alpha_i=1}$ ($\max\emptyset\df 0$).

\begin{lemma}\label{slice-lem}\rm
Let $\BE$ be a measurable vector space and $\Phi:\BE\to \R_+$ be a measurable
function. Let $\ca C$ and $\wt{\ca C}$ be classes of finite random sequences
such that
\begin{equation}\label{slice-1}
\E\Phi \left(x+\sum_{i=1}^m \xi_i x_i \right)\le
\E\Phi \left(x+\sum_{i=1}^m \wt{\xi}_i x_i \right)
\end{equation}
for every integer $m$, $\set{x_i}\subset \E$, $(\xi_i)\in \ca C$, and
$(\wt{\xi_i})\in \wt{\ca C}$.  Let $\BX$ and $\widetilde{\BX}$ be $\ca C$- and
$\widetilde{\ca C}$-sliceable $n\times N$ random matrices, respectively, and
$\bm f=(f_\ba):\ba \subset [1,N])$ be an $\BE$-valued symmetric function. Then
\begin{equation}\label{slice-mult}
\E\Phi \left(\lll\bm f\BX^\otimes    \rrr \right)\le
\E\Phi \left(\lll\bm f\wt{\BX}^\otimes\rrr \right)
\end{equation}
\end{lemma}                                     
\Proof
The statement will be proved by induction with respect to $n$. Without loss of
generality, we may assume that $\BX$ and $\wt{\BX}$ are independent and defined
on a product space, and functions $f_\ba(\ba\bm i)$ vanish unless
$i_1<i_2<i_3\ldots$.

For $n=1$, \refrm{slice-mult} coincides with \refrm{slice-1}. Suppose that
\refrm{slice-mult} holds for every $\ca C$-sliceable matrix $\BX$ and every
$\wt{\ca C}$-sliceable matrix $\wt{\BX}$. We note the decomposition:
\[
f_\ba(\ba\bm i)=
f_\ba (\ba\bm i)\Ib{i_{\ba^*}\le n-1}+
f_\ba(\ba \bm i)\Ib{i_{\ba^*}=n}
\]
Hence, denoting 
\[
f_\ba^{(n)}(\ba\bm i)=
f_\ba(\ba \bm i)\Ib{i_{\ba^*}\le n},
\]
and 
\[
\wt{f}_{\ba\setminus\set{\ba^*}}(\ba\bm i)=
\wt{f}^{(n-1)}_{\ba\setminus \set{\ba^*}} \Ib{i_{\ba^*}=n}
\]
we have 
\begin{eqnarray*}
&\lll \bf \BX^\ba\rrr=
\Sum_{\ba \subset[1,N]} 
\ll f_\ba\BX^{\otimes\ba}\rr&\\
&=\Sum_{\ba \subset[1,N]} 
\ll f_\ba^{(n-1)} \BX^\ba\rr+
\Sum_{\ba \subset[1,N]} 
\ll \wt{f}^{(n-1)}_{\ba\setminus\set{\ba^*}}
\BX^{\otimes \ba\setminus\set{\ba^*}}\rr X_{\ba^*, n}.&\\
\end{eqnarray*}
Now, we use the Fubini's theorem, combined, first, with the inductive
assumption, and then, with condition \refrm{slice-1}. This completes the
proof.\QED

\begin{remarks}\label{slice-rem}\rm
Several special cases and variations of the above lemma will be of particular
interest.
\begin{enumerate}
\item 
Let $\BX=[\BX_1,\ldots,\BX_n]$, where $\BX_k$ is an $n\times k$ random matrix,
and $\wt{\BX}$ have the same structure. Assume that both
matrices are $\ca C$- and $\wt{\ca C}$-sliceable, respectively, and let
$N=1+\ldots +n$. Let $\BE$ and $\Phi$ be as in Lemma \ref{slice-lem}.
\begin{enumerate}
\item
If condition \refrm{slice-1} is fulfilled then, for every symmetric function
$\bm f =(f_k:0\le k\le n)$, we have
\begin{equation}\label{slice-mult-k}
\E\Phi\left(\sum_{k=0}^n \ll f_k \BX_k^\otk \rr \right)\le 
\E\Phi\left(\sum_{k=0}^n \ll f_k \wt{\BX}_k^\otk \rr \right)
\end{equation}
\item
Assume, additionally, that columns of $\BX$ and $\wt{\BX}$ are independent, and
the classes $\ca C$ and $\wt{\ca C}$ are closed under independent extensions
(i.e., if $\bm\xi,\bm \xi'\in\ca C$, and $\bm\xi $ is independent of $\bm
\xi'$, then $(\bm\xi,\bm\xi')\in\ca C$).  Then the following inequality is
sufficient  for \refrm{slice-mult-k}.
\begin{equation}\label{slice-1-1}
\E\Phi(x+\xi y)\le \E\Phi(x+\wt{\xi} y)
\end{equation}
\end{enumerate}
\item 
Let $0<q<p<\infty$ and $\ca C,\wt{\ca C}, \BX,\wt{\BX}$ be as in the lemma or
as  in the special case described above (in Remark \ref{slice-rem}.1). Assume
that
\begin{equation}\label{slice-1-pq}
\|x+\sum_{i=1}^m \xi_i x_i\|_p\le
\|x+\sum_{i=1}^m \wt{\xi}_i x_i\|_q.
\end{equation}
Then
\begin{equation}\label{slice-mult-pq}
\|\lll\bm f\BX^\otimes    \rrr\|_p\le
\|\lll\bm f\wt{\BX}^\otimes\rrr\|_q,
\end{equation}
and, in the special case (Remark \ref{slice-rem}.1),
\begin{equation}\label{slice-mult-k-pq}
\|\sum_{k=0}^n \ll f_k \BX_k^\otk \rr\|_p\le 
\|\sum_{k=0}^n \ll f_k \wt{\BX}_k^\otk \rr\|_q.
\end{equation}
Consider the assumption in Remark \ref{slice-rem}.2. Then
\refrm{slice-mult-k-pq} is fulfilled provided
\begin{equation}\label{slice-1-1-pq}
\|x+\xi y\|_p\le \|x+\wt{\xi} y\|_q
\end{equation}
holds. If there exists a constant $c$ such that $\wt{\xi}=c\xi$, relations
\refrm{slice-1-pq}--\refrm{slice-1-1-pq} are called {\em hypercontraction}
inequalities,  and $\xi$ is called a hypercontractive random variable. Gaussian
and Rademacher random variables are hypercontractive with constants
$c=c_{p,q}=((p-1)/(q-1)^{1/2}$, $1<q<p<\infty$ (cf.
\cite{Bor:chaos,Gro,KraS:hyp,KwaS}).  A symmetric $\alpha$-stable random
variable is hypercontractive in any normed space with exponents $q,p\in
(h_\alpha,\alpha)$, where $h_\alpha=0$, for $\alpha\le 1$, and $h_\alpha<1$,
for every $\alpha<2$ \cite{Szu:hyps}.
\end{enumerate}
\end{remarks}
\subsection{Tail estimates}
In \cite{KwaW:book} the following relation between two $\BE$-valued random
vectors $X$ and $Y$ is called the {\em $\Phi$ domination of $X$ by $Y$}:
\[
\E\Phi(x+X)\le \E\Phi(x+Y),\quad x\in \BE.
\]
In case when $\Phi(\cdot)=\|\cdot\|^p$, for some $p>0$, we will use the phrase 
``$(\BE+L^p(\BE))$-domination (to distinguish the notion from the comparison of
moments).

The fulfillment of $\Phi$-domination yields immediately the same relation for
sums of finite copies of $X$ and $Y$. In \cite{Szu:dechom} we used that fact to
prove that the $\Phi$-domination of two type of random chaoses generated by
hypercontractive random variables implies the tail domination. We will rephrase
that result in a more general manner, pointing out the assumptions needed for
the fulfillment of the tail decoupling.

\begin{theorem}\label{tail}\rm 
Let a class $\ca X$ of random vectors $X$ be $(\BE+L^p(\BE))$-dominated by a
class $\ca Y$ of random vectors $Y$ in $c_0$ (or, equivalently, in every
separable Banach space). Let $\ca Y$ satisfy the, so called,
Marcinkiewicz-Paley-Zygmund (MPZ) condition, i.e.
\[
m=\sup_{Y\in \ca Y}\Frac{\|Y\|_p}{\|Y\|_q}<\infty
\]
for some (equivalently, all) $q<p$. Then, for some constants $c,C>0$, for every
$Y\in \ca Y$, there exists $t_0=t_0(\ca L(\|Y\|)$, such that
\[
\P(\|X\|>ct)\le C\P(\|Y\|>t),\quad t\ge t_0.
\]
If, additionally,  $\ca Y$ is bounded in $L^0(\BE)$, then the class $\ca X$ is
tail-dominated by the class $\ca Y$ (i.e., the number $t_0$ above does not
depend on a particular choice of $Y\in \ca Y$.
\end{theorem}                                         

We omit the proof, since its steps are exactly the same as steps in  the proof
of Theorem 5.3 in \cite{Szu:dechom}. Also, as in \cite{Szu:dechom}, we obtain
immediately the following corollaries.
\begin{cor}\label{cor:tail}\rm
Let assumptions of Theorem \ref{tail} be fulfilled, including the boundedness
of $\ca Y$ in $L^0(\BE)$.
\begin{enumerate}
\item
Let $\varphi:R_+\to\R_+$ be an increasing function of moderate growth, and
$\phi(0)=0$. Then, for some $C'>0$,
\[
\E\varphi(\|X\|)\le C'
\E\varphi(\|Y\|),\qquad X\in \ca X, \quad Y\in \ca Y,
\]
If the growth is not moderate, then we still preserve 
the implication
\[
\E\varphi(\|Y\|)<\infty \imp 
\E\varphi(\|cY\|)<\infty,\qquad X\in \ca X, \quad Y\in \ca Y,
\]
for some universal constant $c>0$.
\item 
The $L^0$-boundedness of $\ca Y$ implies the $L^0$-boundedness of $\ca X$. If
$\ca Y$ is tight, so is $\ca X$.
\item
The domination in the sense of tightness also holds in any separable Fr\'echet
(i.e., metrizable complete locally convex) space, with the topology generated
by a countable family of seminorms (cf. \cite{Rud}), provided the
$(\BE+L^p(\BE))$-domination is fulfilled and the uniform
Marcinkiewicz-Paley-Zygmund condition is fulfilled for all seminorms.
\end{enumerate}
\end{cor}
It is clear how this pattern applies to decoupling inequalities. If a (lower or
upper) decoupling inequality holds in every Banach space, and a random chaos is
induced by hypercontractive random variables, then the same type of decoupling
holds by means described in the above corollary.



\subsection{Lower decoupling}
We assume in this subsection that $\L$ is an Orlicz space $L^\varphi$ such that
$\varphi$ satisfies a strong convexity condition
\begin{equation}\label{Stass-by1}
\mbox{\em for some $a<1$, $\varphi^a$ is convex} 
\end{equation}
Note that \refrm{Stass-by1} means that, for some $p>1$, $\lim_{t\to\infty}
\varphi(t)/t^p=\infty$. In particular, for a moderately increasing $\varphi$
(i.e., for separable $L^\varphi$),  $\L$ is uniformly convex.
We begin with an auxiliary result.
\begin{lemma}\label{l:sup}\rm
Let $\L$ and $\varphi$ satisfy \refrm{Stass-by1}.  Let $\bm \theta=(\theta_i)$
be a sequence of integrable \iid random variables. Put
$S_n=\theta_1+\ldots+\theta_n$ and $\xi=\sup_n|S_n|/n$, and let $(\eps,
\eps_1,\eps_2,\ldots)$ be a Rademacher sequence independent of $(\theta_i)$.
Then, there exists a constant $c_\varphi$, such that
\begin{itemize}
\item[(i)] For every $x,y\in \BE$
\[
\E\varphi(\|x+\eps \xi y\|)\le 
\E\varphi(\|x+c_\varphi\eps \theta y\|).
\]
\item[(ii)] For every $n\in\N,\,x,x_1,\ldots,x_n\in\BE$,
\[
\E\varphi(\|x+\sum_{i=0}^n\eps_i \frac{S_i}{i} x_i\|)\le 
\E\varphi(\|x+c_\varphi\theta\sum_{i=0}^n\eps_i x_i\|).
\]
\end{itemize}
\end{lemma}
\Proof
Assertion (ii) follows immediately from (i), the Fubini's theorem, and  the
contraction principle for a Rademacher sequence.

We will prove (i). The function $[0,\infty)\ni t\mapsto
\psi(t)=\E\varphi(\|x+\eps t y\|)-\varphi(\|x\|)$ is convex and increasing.
Hence
\[
\E\psi(\xi)\le C\E\psi(\theta),
\]
since the sequence $M_1=S_n/n,M_2=S_{n-1}/(n-1),\ldots,
M_{n-1}=S_2/2,M_n=S_1=\theta_1$ forms a martingale with respect to the natural
filtration.

It is an elementary exercise to prove that the following transformations
inherit property \refrm{Stass-by1}: the shift $\phi-a$, the composition
$\phi\circ\psi$ with another convex function, averages $\int
\phi_\omega(\cdot)\mu(d\omega)$ with respect to probability measures $\mu$ and
a (measurable) family $\{\phi_\omega\}$ of functions with property
\refrm{Stass-by1}.  Hence the function 
$[0,\infty)\ni t\mapsto \E\phi(\|x+\eps t y\|)-\phi(\|x\|)$ has property
\refrm{Stass-by1}, where $x,y\in \BE$ and $\eps$ is a Rademacher random
variable.
Therefore, by Doob's inequality,
\begin{equation}\label{Doob}
\P(\xi>t)=\P(\psi(\xi)>\psi(t))\le 
\frac{\E[\psi^a(\theta);\xi>t)}{\psi^a(t)}.
\end{equation}
Then, we infer from \refrm{Doob} and H\"older's inequality that
\[
\begin{array}{rl}
\displaystyle
\E\phi(\|x+\eps y\xi\|)-\phi(1)
&=\E\psi(\xi)=\Int_0^\infty\P(\xi>t)\,d\psi(t)\\
\displaystyle
\le&\Int_0^\infty\Frac{\E[\psi^a(\theta);\xi>t)}{\psi^a(t)}
\,d\psi(t)\\
\displaystyle
=&(1-a)^{-1}\E[\psi^a(\theta)\psi^{1-a}(\xi)]\\
\displaystyle
\le&(1-a)^{-1}(\E[\psi(\theta)])^a(\E[\psi(\xi)])^{1-a}.\\
\end{array}
\]
Define $a_0=\inf\set{a\in(0,1):\phi^a\mbox{ is convex}}.$
Then, letting $a\to a_0$, and using the convexity, we obtain
\[
\E\psi(\xi)\le (1-a_0)^{-1/a_0}\E\psi(\theta)\le 
\E\psi((1-a_0)^{-1/a_0}\theta).
\]
The lemma has been
proved.\QED
\begin{remark}\label{constant}\rm
The  constant $c=c_{\varphi}$ depends on the exponent $a$, appearing in
\refrm{Stass-by1}, or more precisely, on  $a_0$. Also, $c=\infty$,  if $a_0=1$,
in general.  The convexity assumption concerning  $\varphi$ is necessary for
(i), if we do not restrict the class of distributions of $\theta$. Consider,
for example, $\L=L^1$. Then (i) implies that $\E\sup_i|\theta_i|/i<\infty$, if
$\E|\theta|<\infty$. A symmetric random variable $\theta$ with the tail
$\P(|\theta|>t)=(t\log^2t)^{-1}$, $t\ge e$, produces a quick counterexample.
\end{remark}

We will let the generality of the proof of the following theorem slightly
exceed our current needs. The reason will be explained in Section 5. Recall 
(see \refrm{realdec})
that in the real case the lower decoupling for Gaussian or Rademacher
chaoses holds with a constant $c_k=1/\sqrt{k!}$.

\begin{theorem}\label{dom-by-chaos}\rm
Let the matrix $[\bx,\BX]$ have i.i.d. columns. Let $\varphi$
satisfy \refrm{Stass-by1}. Then the sign-randomized weak lower decoupling
inequality holds, i.e., for  Walsh function $\bm w=(w_k)$, independent of
$[\bx,\BX]$, we have
\begin{equation}\label{e:dom-by-chaos}
\|\Sum_{k=0}^n w_k\ll f_k \BX^{\otk} \rr \|\le 
\|\Sum_{k=0}^n w_k\Frac{(2ck)^k}{k!}
\ll f_k \bx^{\otk}\rr\|, 
\end{equation}
for every symmetric function $\bm f=(f_k)$ vanishing on diagonals, where $c
=c_\varphi$ depends only on the function $\varphi$. If the underlying random
variables are symmetric, then the Walsh functions can be omitted.
\end{theorem}
\Proof 
We begin with the Mazur-Orlicz
polarization formula.
\[
\begin{array}{rl}
\|\San\ll  f_\ba\BX^{\otimes\ba} \rr \|
=&\|\displaystyle\San\ll  f_\ba\rec{|\ba|!}\Sum_{\bb\le \ba}
(-1)^{|\ba-\bb|}\l(\bb\BX\r)^{\otimes\ba} \rr \|\\
=&\displaystyle \|2^{-n}\Sum_{\bb\le [1,n]}
\ll \San \rec{|\ba|!}(-1)^{|\ba-\bb\ba|}  f_\ba
\l(2\bb\ba\BX\r)^{\otimes\ba} \rr \|\\
\le&\displaystyle 2^{-n}\Sum_{\bb\le [1,n]}
\|\San\rec{|\ba|!}(-1)^{|\ba-\bb\ba|}\ll  f_\ba\l(2\bb\ba\BX\r)^{\otimes\ba} 
\rr \|\\
\end{array}
\]
For a fixed $\bb$, we have (cf. \refrm{bax})
\[
\l(\bb\ba\BX\r)^{\otimes\ba}=
\E[\,\l(\ba\BX\r)^{\otimes\ba} \,|\,\BX^\bb \,].
\]
Recall that $w_\ba(\bb)=(-1)^{\ba\bb}$ are Walsh functions. Since the mapping
$\bb\mapsto \bb'=1-\bb$ is measure preserving, hence, by the contractivity of
conditional expectations and Fubini's theorem, we have
\[
\|\San\ll  f_\ba\BX^{(\otimes\ba)} \rr \|\le
\|\San w_\ba \Frac{2^{|\ba|}}{|\ba|!}
\ll  f_\ba\l(\ba\BX\r)^{\otimes\ba} \rr\|.
\]
At this moment we give up the generality and notice that $f\ba$ ($|\ba|=k$)
vanish unless $\ba=[1,k]$.

Now it suffices to apply the Slicing Lemma \ref{slice-lem}.  Let $\bm g=(\bm
g_k)$ be a symmetric function, $\bm g_k:\N^k\to\BE$, \refrm{Stass-by1}  be
fulfilled and $\bm w=(w_k)$ be a Walsh sequence independent of $\BX$. Then, for
a constant $c=c_\varphi$,
\begin{equation}\label{e:slicing}
\|\sum_{k=0}^nw_k\ll \bm g_k
\left(\frac{\bx_1+\ldots+\bx_k}{k}\right)^{\otk} \rr \|
\le\|\sum_{k=0}^n w_k\ll  \bm g_k(c\bx)^{\otk} \rr \|.
\end{equation}

The proof is completed.\QED
\begin{cor}\label{corLD}\rm
Let assumptions of Theorem \ref{dom-by-chaos} be fulfilled, where
$\varphi(t)=t^p,\,p>1$. Denote by $Q(f)$ the coupled, and by $\overline{Q}(f)$,
the decoupled chaos, as appear, respectively, in the right and left hand side
of inequality \refrm{e:dom-by-chaos}. Assume that components $X_i$ of $\bx$ are
hypercontractive, with hypercontractivity constants uniformly bounded away from
0. Then the following conditions are fulfilled.
\begin{enumerate}
\item
There is a constant $C$, depending only on the hypercontractivity constants,
and a sequential constant $\bm c=(c_k)$, depending only on $p$, such that, for
any non-decreasing moderately growing function $\varphi:\R_+\to\R_+$, every
$\bm f\in \ca F_S$,
\[
\E\varphi(\|\ov{Q}(f)\|)\le C 
\E\varphi(\|{Q}(cf)\|) 
\]
(if $\varphi$ does not grow moderately, the finitness of the Orlicz modular is
preserved).
\item
The stochastic boundedness of a family of coupled polynomial chaoses
$\set{Q(f_a):a\in A}$ implies the same for $\set{\ov{Q}(\bm d f_a):a\in A}$,
where $\bm d=\bm c^{-1}$, and $\bm c$ appears in the preceding statement. By
the same token, the tightness of the first family yields the tightness of the
second family.
\end{enumerate}
\end{cor}
\Proof
It suffices to interpret appropriately Corollary \ref{comp}. \QED

\subsection{Reduction to Rademacher chaoses}
We will focus on a search of reasonably wide classes of Banach spaces,
which support the exponential upper decoupling.
Recall that $\bx\in \ca{UD}=\ca{ UD} (\BE;\Phi;\ca F)$ ($\bx$ satisfies the
upper decoupling inequality), if
\[
\E\Phi(\lll \bm f \bx^\otimes \rrr)\le
\E\Phi(\lll \bm f (c \BX^\otimes \rrr)
\]
for every function from class $\ca F$.  The most important are classes ${\ca
F}_S$, of  symmetric functions,  and ${\ca F}_T$, of tetrahedral functions.
Denote by $\mu=\ca L(\bx)$ the distribution of a sequence $\bx$. Let
$\varphi:\R_+\to\R_+$ be a measurable function. Denote by 
$\ca R_U=\ca R_U(\mu,\varphi)$ (respectively, $\ca R_L=\ca R_L(\mu,\varphi)$
the class of Banach spaces such that, for some constant $c>0$, the inequality
\begin{equation}\label{Usumy}
\E\varphi(\|x+X\sum_i\eps_i x_i\|)\le 
\E\varphi(\|x+c\sum X_i x_i\|),
\end{equation}
(respectively,
\begin{equation}\label{Lsumy}
\E\varphi(\|x+\sum X_i x_i\|)\le 
\E\varphi(\|x+c X\sum_i\eps_i x_i\|)
\quad )
\end{equation}
is fulfilled, for every $x\in \BE$, and for all finite families
$\set{x_i}\subset \BE$.

The following result shows the importance of the introduced classes. Its proof
is  a direct consequence of the Slicing Lemma \ref{slice-lem}, and, for
tetrahedral functions, of  equality \refrm{Vcontr} from Proposition
\ref{condextensor}.
\begin{proposition}\label{dec-up}\rm 
Let $\bx $ be a sequence of independent symmetric random variables, and $\BX$
be a matrix whose columns are independent copies of $\bx$, $\bm f=(f_\ba)$ be a
symmetric function with values in $\BE$,  $\varphi:\BE\to \R_+$ be a measurable
function.  Denote by $\beps=(\eps_i)$ be a Rademacher sequence independent of
$\bx$, and by $\Bbb S$,  a Rademacher matrix, independent of $[\bx,\BX]$.  If
$\BE\in \ca R_U(\mu,\varphi)$ (respectively, $\BE\in \ca R_L(\mu,\varphi)$),
then
\[
\E\varphi \left(\left\|\lll \bm f (\bx {\Bbb S})^\otimes\rrr \right\|\right)\le
\E\varphi \left(\left\|\lll \bm f (c\BX)^\otimes\rrr \right\|\right)
=\E\varphi \left(\left\|\lll \bm f (c{\Bbb S} \BX)^\otimes\rrr \right\|\right)
\]
(respectively, the converse implication is valid, with $c$ replaced by
$c^{-1}$).  If, additionally, $\varphi$ is convex, the latter inequality is
fulfilled also for tetrahedral functions (respectively, the fulfillment of the
latter inequality for tetrahedral functions implies the same, for symmetric
functions).
\end{proposition}                                         

\subsection{Limitations of the reduction}
Inequalities \refrm{Lsumy} and \refrm{Usumy} may fail in some Banach spaces,
and for some random sequences.  First, we note the following immediate
consequence of Proposition \ref{dec-up}.  %
\begin{lemma}\label{dec-up easy}\rm 
Let $\bx,\BX$, and $\varphi$ be as in Proposition \ref{dec-up}. Let $\BE\in \ca
R_U(\mu,\varphi)$ (respectively, $\BE\in \ca R_L(\mu,\varphi)$. Then
\begin{equation}\label{e:Ueasy}
\E\varphi(\|x+\sum_{j=1}^m\sum_{i=1}^n X_j \eps_{ij} x_{ij}\|)\le 
\E\varphi(\|x+c\sum_{j=1}^m\sum_{i=1}^n X_{ij} x_{ij}\|)
\end{equation}
(respectively,
\begin{equation}\label{e:Leasy}
\E\varphi(\|x+\sum_{j=1}^m\sum_{i=1}^n X_{ij} x_{ij}\|)
\le 
\E\varphi(\|x+c\sum_{j=1}^m\sum_{i=1}^n X_j \eps_{ij} x_{ij}\|)
\quad)
\end{equation}
for every $m,n\in\N$, and every matrix $[x_{ij}]$ of vectors of $\BE$.
\end{lemma}                                         

\begin{proposition}\label{no R}\rm
Let $\bx,\BX$ be as in Proposition \ref{dec-up}, $\varphi(t)=t^2$, and
$\bm G=(G,G_1,G_2,\ldots)$, and $\bm Y=(Y,Y_1,Y_2,\ldots)$, respectively, be a
sequence of i.i.d. standard normal, and a sequence of i.i.d. exponential random
variables with parameter $2$, respectively, such that $\bm G$ and $\bm Y$ are
independent of $\bx$.
\begin{enumerate}
\item
Let $X$ be 
square integrable. Assume that $\BE\in \ca R_U(\mu,\varphi)$. Then the
following conditions are fulfilled.
\begin{itemize}
\item[{\rm (i)}] 
\begin{equation}\label{twopoint-theta-gamma}
\E\|x+X G y\|^2
\le \E\|x+c G y\|^2,\qquad x,y\in \BE;
\end{equation}
\item[{\rm (ii)}] 
\begin{equation}\label{twopoint-YG}
\E\|x+\eps Yy\|^2 \le \E\|x+c'G y\|^2, \qquad x,y\in \BE,
\end{equation}
where $c'$ may be a new constant;
\item[{\rm (iii)}] 
\begin{equation}\label{twopoint-YG-sums}
\E\|x+\sum_i\eps_i Y_i x_i\|^2 \le 
\E\|x+c'\sum_i G_i x_i\|^2, \qquad x,y\in \BE
\end{equation}
\item[{\rm (iv)}] 
$\BE$ does not contain isomorphic copies of $\ell^\infty_n$, uniform in $n$,
neither it contains two-dimensional subspaces isometric to $\ell^\infty_2$ or
$\ell^1_2$.
\end{itemize}
\item
Assume that $\BE\in \ca R_L(\mu,\varphi)$, where 
$\varphi$ be a nondegenerate nondecreasing function such that
$\limsup_{t\to\infty} \varphi(t) e^{-at^2}=0$ for some $a>0$.   Then $X$ is 
square integrable. 
\end{enumerate}
\end{proposition}                                         
\Proof
1. Inequality \refrm{twopoint-theta-gamma} follows by Lemma \ref{dec-up easy}
and the (real) Central Limit Theorem.  The passage to the limit can be
justified by a routine uniform integrability argument (cf., e.g.,
\cite[Theorems 5.3 and 5.4]{Bil}).  By the same token, one may assume that $X$
in  \refrm{twopoint-theta-gamma} has the normal distribution.  An exponential
random variable $Y$ with intensity $\lambda=2$ has the tail comparable to the
tail of the product of two independent Gaussian random variable (cf. \cite[pp.\
243-244]{Yos}), which proves \refrm{twopoint-YG}, in view of Corollary
\ref{comp}.  Estimate \refrm{twopoint-YG-sums} follows by the Fubini's theorem
and iteration.

Suppose $\BE=\ell^\infty_2$ (i.e., $\BE$ is just $\R^2$ with sup-norm). Take
orthogonal $y,x$ in
\refrm{twopoint-YG} with $\|x\|=1, \|y\|=1/u<1$, and subtract 1 from both sides
of \refrm{twopoint-theta-gamma}. Then the right hand side is of order
$\ex{-u^2/2}$, while the left hand side is of order $\ex{-2u}$, for $u\to
\infty$, which produces  a contradiction.

Inequality \refrm{twopoint-YG-sums} yields the domination of sums of
symmetrized independent exponential random variables by sums of independent
Gaussian random variables. Clearly, this is impossible in $\ell^\infty_n$ (it
suffices to take orthogonal $x_i$'s, and apply the classical estimates for
suprema of independent random variable, cf., e.g., \cite[Lemma V.3.2]{VCT}).

Finally, inequality \refrm{twopoint-YG} does not hold in $\ell^1_2$, since by
the Ferguson-Hertz embedding theorem every two-dimensional normed space can be
isometrically embedded into $L^1$ (cf. \cite{Fer,Her}, see also \cite{KwaS}).
One can construct a direct counterexample, too.
\vspace{5pt}

2. By choosing
$x_i=t_n x/\sqrt{n}$, $i=1,\ldots,n$, where $\|x\|=1$, in the defining 
inequality of the class $\ca R_L$, we infer that, for some constant $c'$
\begin{equation}\label{LD to Rad in L2 only}
\E\varphi(t_n \frac{|\sum_{i=1}^n X_i|}{\sqrt{n}})\le 
c' \E\varphi(t_n \frac{|\sum_{i=1}^n \eps_i|}{\sqrt{n}}),
\end{equation}
for every real sequence $t_n\to 0$. Because of the regular variability of
$\varphi$ at $\infty$, the right hand side converges to $0$, hence the sequence
$(\sum_{i=1}^n X_i/\sqrt{n})$ is bounded in $L^0$ (i.e., tight), by Chebyshev's
inequality. This is possible only if $X\in L^2$.
\QED
\subsection{Reduction in some spaces}
\subsubsection{\em Rademacher versus Gaussian chaoses}
So far, we have established a class of Banach spaces, where the exponential
upper decoupling inequality is fulfilled. Now, we will show that class is
reasonably wide. In general, the upper decoupling may depend also on
distributions of involved  random variables (whether it does, is an open
question at this time). Before we proceed further, in order to avoid
unnecessary redundancy, we will determine some dependence (far from being
complete) between decoupling inequalities for random chaoses spanned by random
variables with different distributions.
\begin{proposition}\label{closeness}\rm
\begin{itemize}
\item[{\rm (i)}] The class $\ca{ UD}_S=\ca{ UD}  (\BE;\|\cdot\|^p,\ca F_S)$ is
closed under products and sums of i.i.d. sequences, i.e., if $\bx $ and  
$\bx'$ are equidistributed independent sequences, and $\bx, \bx'\in\ca{ UD}_S$
with a constant $c$, then $\bx\bx'=(X_iX'_i)\in \ca {UD}$ and $\bx+\bx'\in \ca
{UD}_S$ with the constant $c$.
\item[{\rm (ii)}] 
In addition to the above properties,  the class $\ca {UD}_T$  is also closed
under linear combinations of independent sequences, i.e., if $\bx $ and  $\bx'
$ are independent sequences, and $\bx, \bx'\in\ca{ UD}_T$ with constants
$c,c'$, then, for every numerical sequences $\bm a$ and $\bm b$, $\bm a\bx+\bm
b \bx'\in \ca {UD}_T$ with the constant $c$.
\item[{\rm (iii)}]
Denote $\psi(t)=\psi_{x,y}(t)=\E\Phi(x+\eps t y)$. If $\bx^{(m)}\in\ca{UD}$
with the same constant $c$, the distributions of $\bx^{(m)}$ converge weakly to
the distribution of $\bx$, and the family $\set{\psi(\bx^{(m)})}$ is uniformly
integrable, then $\bx\in\ca{UD}$ with a constant which is less or equal $c$.
\item[{\rm (iv)}]
If $\Phi$ is convex, then ${\cal UD}_S\subset{\cal UD}_T$
\end{itemize}
\end{proposition}                                         
\Proof
The closeness under the product is easy to see and follows immediately by
Fubini's theorem.

In order to prove the additivity in (i), let us consider a $(2n\times
2n)$-matrix
\[
\left[
\begin{array}{cc}
\BX &\BY\\
\BX'&\BY'\\
\end{array} 
\right] 
\]
where $\BY$ and $\BY'$ are independent copies of $\BX$, and the sequence
$(\bx,\bm Y)$, where $\bm Y$ is an independent copy of $\bx$. Then it suffices
to change the enumeration of arguments of functions $f_k(\cdot)$, putting, in
particular $\bm f_k=0$ for $k\in [n+1,2n]$.

For additivity in (ii), we rather use the following $(2n^2\times n)$-matrix
\[
\left[
\begin{array}{ccccc}
\bx_1 &*     &*     &\ldots&*     \\
\bx'_1&*     &*     &\ldots&*     \\
   {*}&\bx_2 &      &\ldots&*     \\ {*}&\bx'_2&      &\ldots&*     \\ {*}&*
&\bx_3 &\ldots&*     \\ {*}&*     &\bx'_3&\ldots&*     \\
\ldots&\ldots&\ldots&\ldots&\ldots\\
   {*}&*     &*     &\ldots&\bx_n \\ {*}&*     &*     &\ldots&\bx_n'\\
\end{array}\right],
\]
where the symbols $*$ indicate the presence of mutually independent copies of
corresponding portions of columns.

Note that the lack of symmetry assumption in (ii) enables us to use arbitrary
sequential multipliers $\bm a$ and $\bm b$, while both sequences must be
constant under the symmetry assumption.

Assertion (iii) follows from basic properties of weak convergence (cf.
\cite[Theorems 5.3 and 5.4]{Bil}).

Assertion (iv) follows from \refrm{Vcontr}.
\QED

\begin{cor}\label{RadtoGauss}\rm 
If the upper decoupling inequality for symmetric  (or triangular functions) is
satisfied for some zero-mean probability law with finite variance, e.g., by a
Rademacher random variable, then it is satisfied by the Gaussian law.
\end{cor}
\Proof 
The assertion follows from Proposition \ref{closeness}, (i) or (ii), the
Central Limit Theorem, and Proposition \ref{closeness}(iii).
\begin{cor}\label{Radtoexp}\rm
Let the assumptions of Theorem \ref{RadtoGauss} be fulfilled. Then the upper
decoupling inequality of the same type (i.e., either for symmetric or
triangular functions) is satisfied by all symmetrized 
${\rm Gamma}(m)$-distributions, $m=1$ (exponential law), $2,3,\ldots$
\end{cor}
\Proof 
Indeed, the product of two independent Gaussian random variables is comparable
to a random variable with exponential distribution (cf. e.g. 
\cite[pp.243-244]{Yos}). Hence the assertion follows by Proposition
\ref{closeness}.\QED

\subsubsection{\em Banach lattices}
Let $\BE$ be a Banach lattice. Then, for every continuous positive homogeneous
function $\psi:\R^n\to\R$, the expression $\psi(x_1,\ldots,x_n)\in\BE$,
$x_1,\ldots,x_n\in \BE$, is well defined, in particular,
\[
(\sum_{i=1}^n |x_i|^p)^{1/p};\quad (\E|\sum_{i=1}^n x_i\theta_i|^p)^{1/p},
\]
where $\theta_i\in L^p$ are real random variables, and $0<p\le\infty$
(\cite{Kri}, also see \cite{LinT}). Any inequality or equality that is valid in
the real case, carries over to Banach lattices (when $\BE$ is a space of
functions, these constructions, in general, can be viewed
pointwise, both intuitively and rigorously).

Recall the Krivine's notion of type $\ge p$ and $\le p$ ($p$-convex and
$p$-concave in \cite{LinT}). A Banach lattice is said to be of type $\ge p$
(respectively, of type $\le p$), $1\le p \le \infty$ if
\[
\|(\sum_{i=1}^n |x_i|^p)^{1/p}\|\le C
(\sum_{i=1}^n \|x_i\|^p)^{1/p}
\]
(respectively, the inverse inequality holds). These properties refer to a
degree of convexity of the unit sphere, compared to the unit sphere  in $L^p$.
For example, $L^r$ is of type $\le p$, when $r\le p$, and of type $\ge p$, when
$r\ge p$, $0<r\le \infty$.

\begin{theorem}\label{UDlust}\rm
Every Banach lattice $\BE$ of type $\le q<\infty$ and type $p>1$ admits an
equivalent norm for which both upper and lower (and both symmetric and
tetrahedral)  decoupling inequalities hold for the Rademacher (hence Gaussian)
law, by means of comparison in $L^r$, $1<r<\infty$. More precisely, for such a
norm, there exists a constant $c$ such that
\[
\|\sum_{k\ge 0} \ll Q_k(\bm f/c)\| \le
\|\sum_{k\ge 0} \ll \overline{Q_k}(\bm f)\| \le
\|\sum_{k\ge 0} \ll Q_k(c \bm f)\|,
\]
where
\[
\overline{Q}_k(f_k)=\left\{
\begin{array}{rl}
\Rec{\sqrt{k!}}\ll f_k \BX^\otk\rr& \mbox{if $f_k$ is symmetric}\\
\ll f_k \BX^\otk\rr& \mbox{if $f_k$ is tetrahedral},\\
\end{array}
\right.
\]
and $Q_k(f_k)=\ll f_k \bx^\otk\rr$ in both cases.
\end{theorem}                                         
\Proof
Essentially, we reduce the problem to the situation on the real line (cf.
\refrm{realdec}).  By Figiel and Johnson  theorem (\cite{FigJ}, see also
\cite[Theorem 1.d.8]{LinT}), a Banach lattice $\BE$, which is of type $\ge p$
and $\le q$, $1<p\le  q< \infty$, admits an equivalent norm, making both
constants $C$, appearing in the definition, equal to 1. So, assume that is the
case.  Then, we have
\[
\|(\E|\sum_{i=1}^n x_i\theta_i|^p)^{1/p}\|\le
(\E\|\sum_{i=1}^n x_i\theta_i\|^p)^{1/p}
\]
and
\[
(\E\|\sum_{i=1}^n x_i\theta_i\|^q)^{1/q}\le
\|(\E|\sum_{i=1}^n x_i\theta_i|^q)^{1/q}\|
\]
for any collection of suitably integrable random variables $(\theta_i)$. We
will apply both inequalities to Rademacher chaoses and use the
hypercontractivity of Rademacher chaos. Denote 
$c_{r,q}=\max(1, ((r-1)/(q-1))^{1/2})$. The proof is similar in the symmetric
and tetrahedral case, and also for the upper and lower decoupling. We will give
details only in one case, say, for tetrahedral $\bm f$ and the upper
inequality. We have
\[
\begin{array}{rl}
\left(\E \left\|\lll \bm f \beps^\otimes
\rrr \right\|^r \right)^{1/r}
\le &\left\|\left(\E|\lll\bm f
(c_{r,q}\beps)^\otimes\rrr|^q\right)^{1/q}\right\|\\ 
\le \left\|(\E|\lll\bm f
(c_{r,q}c_{q,2}\beps)^\otimes\rrr|^2)^{1/2}\right\| 
= &\left\|\lll|\bm f|^2
(c_{r,q}c_{q,2})^\otimes\rr)^{1/2}\right\|\\
\le \left\|(\E|\lll\bm f 
(c_{r,q}c_{q,2}\Bbb S)^\otimes\rrr|^2)^{1/2}\right\|
\le &\left\|(\E|\lll\bm f (c_{r,q}c_{q,2}c_{2,p}
\Bbb S)^\otimes\rrr|^p)^{1/p}\right\|\\
\le \left(\E \left\|\lll \bm f (c_{r,q}c_{q,2}c_{2,p}
\Bbb S)^\otimes\rrr \right\|^p \right)^{1/p} 
\le& \left(\E \left\|\lll\bm f (c_{r,q}c_{q,2}c_{2,p}c_{p,r}
\Bbb S)^\otimes\rrr \right\|^r \right)^{1/r}.\\ 
\end{array}
\]
Other cases follow by an almost verbatim argument.\QED

The class of Banach lattices, appearing in the theorem,  can be enlarged to
uniformly convex spaces with a local unconditional structure (LUST) (i.e., such
that $\BE$ can be embedded into the dual of a Banach lattice,  cf. , e.g.,
\cite{GorL}). In fact, the context of Banach lattices, as appear in the
theorem, makes the problem of decoupling rather trivial. Any tetrahedral
(respectively, symmetric and of finite order) Rademacher or Gaussian chaos is
exponentially equivalent, in the sense of the introduced domination in any
$L^r$, $1<r<\infty$, to an infinite (respectively, finite) Rademacher or
Gaussian sum
\[
\sum_{k=0} \sum_{\bm i\in \N^k} f_k(\bm i) X_{\bm i},
\]
where $\set{X_{\bm i}:\bm i\subset \N^\N, \bm i\mbox{ finite}}$ is a family of
independent Rademacher or Gaussian variables. In particular, for the
aforementioned class of Banach spaces, after a renorming, in a trivial manner
an infinite order contraction principle holds for Rademacher or Gaussian
chaoses
\[
\E\|\lll \bm{fg}\BX^\otimes\rrr\|^r\le 
\E\|\lll \bm f(c\BX)^\otimes\rrr\|^r,
\]
where $c$ is a suitable constant, and $\bm g=(g_k)$, $g_k:\N^k\to [-1,1]$ is an
arbitrary measurable function.  Such the contraction principle fails in
general, e.g., if $\BE=c_0$, and even for a single $k$-homogeneous component,
$k\ge 2$ \cite{KraS:sum}.
\begin{remark}\rm
A related procedure can be applied for random variables with sufficiently high
moments. That is, if $\BE$ is as in Proposition  \ref{UDlust}, and
$\theta\in L^{r}$, $r\ge 2$,  is a symmetric random variable such that
$r>q_0=\inf\set{q:\BE \mbox{ is of type } \le q}$, then $\theta$ is
hypercontractive with constants $c_{r,q}(\theta)=
\|\theta\|_r/\|\theta\|_q c_{r,q}$  \cite{KraS:hyp}, which would replace
constants $c_{r,q}$ in the proof. By a similar argument to the one used in the
proof of Proposition \ref{no R}, one can show that $L^q\notin \R_U(\ca
L(X),L^r)$, if $q\ge q_0>r_0\df\sup\set{r:\theta\in L^r}$ ($L^q$ in the latter
formula can be replaced by any Banach spaces containing isomorphic copies of
finite dimensional spaces $\ell^q_n$, uniform in $n$). This may suggest that
the upper decoupling inequality fails in such spaces yet the problem remains
open.
\end{remark}

We do not know whether the upper decoupling inequalities for Gaussian and
Rademacher chaoses are equivalent. For tetrahedral functions, even in the
homogeneous case, the lower decoupling inequality may fail (cf. Bourgain's
example included in \cite{McCT:bandec}, or \cite[Section 6.9]{KwaW:book}).
\section{STABLE CHAOSES}
\subsection{Auxiliary definitions and inequalities}
In this section $X$ denotes a symmetric standard $\alpha$-stable ($S\alpha
S$, in short) random variable, i.e., $\E\ex{itX}=\ex{-|t|^\alpha}$, and $Y$
denotes a symmetric $\alpha$-Pareto random variable, i.e.,
$\P(|Y|>t)=t^{-\alpha},\,t\ge 1$ ($S\alpha P$, in short) . It is known that
tails of $S\alpha S$ and $S\alpha P$ random variables are comparable, i.e.
$\P(|X|>t)\le K \P(|Y|>Kt)$ and  $\P(|Y|>t)\le K \P(|X|>Kt)$, $t\ge t_0$
(we may assume that the above estimate are valid for all $t\ge 0$). Hence,
for  tetrahedral or symmetric infinite order polynomial $Q(\bm f,\cdot)$ we
have
\begin{equation}\label{tails}
\|Q(\bm f, Y/c)\|_p\le 
\|Q(\bm f, X)\|_p\le 
\|Q(\bm f, c Y)\|_p. 
\end{equation}
This remark will enable us to switch freely (in the sense of the exponential
domination) between stable and Pareto chaoses, and benefit from algebraic
properties of stable random variables, or analytic properties of Pareto random
variables.
We will need also the following estimate (cf.
\cite{KraS:hyp,Szu:hyps,Szu:rpconvex} for similar inequalities. 

\begin{lemma}\label{estimate}\rm
Let $0<s<\alpha<2$. There exists a constant $a=a(\alpha,s)$ such that, for
every sequence of i.i.d. $S\alpha S$ (or $S\alpha P$) random variables, the
inequality
\[
(\|x\|^\alpha+a\sum_i\|x_i\|^\alpha)^{1/\alpha} \le
\|x+\sum_i X_i x_i\|_s
\]
is fulfilled, for all $x,x_1,x_2,\ldots\in\BE$.
\end{lemma}                                         
\Proof
We will apply a hypercontractive iteration for Pareto random variables, and
then use the fact that $S\alpha P$ law belongs to the normal domain of
attraction of the $S\alpha S$ law (cf. the aforementioned papers for details).
It suffices to verify the inequality
\begin{equation}\label{Shypdown}
(1+a t^\alpha)^{s/\alpha}\le \E\|x+ Y t y\|^s, 
\end{equation}
where $\|x\|=\|y\|=1, 0<t\le 1$. The inequality follows by combining the
estimate
\[
\E\|x+ Y t y\|^s-1\ge t^\alpha \inf_{\|x\|=\|y\|=1}
\E[\|x+Y y\|^s-1;|Y|\ge 2] \ge t^\alpha \E(|Y|-1|^s-1)_+ 
\]
with the inequality $(1+ t^\alpha)^{s/\alpha}-1\le s/\alpha t^\alpha$, which
holds for all $t\ge 0$. Put
\[
a=\alpha\E(|\,|Y|-1\,|^s-1)_+/s.
\]
This completes the proof.\QED

We will see that the fulfillment of decoupling inequalities may depend on the
convexity and smoothness of the norm.
A norm of a Banach space $\BE$ is called {\em $p$-smooth} (cf.
\cite{Ass,LinT}), $1<p\le 2$, if

\[
(\E\|x+\eps t y\|^2)^{1/2}\le (1 + C t^p)^{1/p}
\]
where $\|x\|=\|y\|=1$, $t>0$ (it suffices to consider only $t\le 1$), and
$\eps$ is a Rademacher random variable. By hypercontractivity, the $L^2$-norm
on the left hand side can be replaced by any $L^s$-norm, $1<s<\infty$.
A Banach $\BE$ is called {\em $p$-smoothable}, if it admits an equivalent
$p$-smooth norm. It will be convenient to  extend trivially the notion of
smoothness to the case $p=1$ (every norm is $1$-smooth).

For  a Banach lattice $\BE$, let
\[
k^0=\inf\set{q:\E\mbox{is of Krivine's type }\le q}\le\infty,
k_0=\sup\set{p:\E\mbox{is of Krivine's type }\ge p}\ge 1.
\]
Clearly, $k_0\le k^0$.
Say that a Banach space is of infinite cotype, if it contains subspaces
isomorphic to $\ell^\infty_n$ uniform in $n$. Otherwise, $\BE$ is said to be a
space of finite cotype. A Banach lattice is of finite cotype if and
only if it is of Krivine's type $\le q$, for some $q<\infty$ \cite{LinT}.

\subsection{Symmetric decoupling}
Let us extract a suitable fragment from Theorem \ref{dom-by-chaos}.

\begin{theorem}[Lower Symmetric Decoupling]\label{LSD}\rm
Let $1<p<\alpha<2$. The lower decoupling inequality in $L^p$ for symmetric
$S\alpha S$ and $S\alpha P$ chaoses is fulfilled  with  constants
$c_k=d^k$, for some $d>0$.
\end{theorem}
The obtained constant is the best we
know, even in the real case. However, for a single homogeneous chaos, the
estimate can be significantly improved (cf., e.g., \cite{DeA}).
Like before, in the Rademacher or Gaussian case, we can prove the upper
decoupling inequality only in some Banach spaces. Surprisingly, the upper
decoupling inequality for nonintegrable stable chaoses is a trivial consequence
of the slicing techniques, and holds in an arbitrary Banach space.  Recall that
any $S\alpha S$ (or $S\alpha P$) random variable is hypercontractive in any
normed space with exponents $q,p\in (h_\alpha,\alpha)$, and $h_\alpha=0$, for
$\alpha\le 1$.
\begin{theorem}[Upper Symmetric Decoupling]\label{USDLat}{\rm}
Let $\BE$ be a Banach space. Consider $S\alpha S$ (or $S\alpha P$) chaoses in
symmetric functions and the norm $L^s$, $h_\alpha<s<\alpha$.
\begin{itemize}
\item[{\rm (i)}]
For $\alpha\le 1$, the exponential upper decoupling inequality for symmetric
$S\alpha S$ (or $S\alpha P$) chaoses holds in every Banach space.
\item[{\rm (ii)}]
Let $\BE$ be $p$-smoothable, $1\le p\le 2$, and $0<\alpha\le p$. Then $\BE\in
\ca R_U(\mu,s)$. In particular, if an upper symmetric decoupling inequality
holds for Rademacher chaoses, with constants $(c^{(R)}_k)$, then, for an
equivalent norm, an upper decoupling inequality for symmetric $S\alpha S$ (or
$S\alpha P$) chaoses holds, with constants $(a c^{(R)}_k)$, where
$a=a(\alpha,s,\|\cdot\|_\BE,p)$.
\item[{\rm (iii)}]
Let $\BE$ be a Banach lattice, $0<s<\alpha<2$.  If $k^0<\infty$, and
$\alpha>k^0$ or $k_0>\alpha$, then $\BE\in \ca R_U(\mu,s)$, and the upper
decoupling inequality with exponential constants $c_k=a^k$ holds.
\end{itemize}
\end{theorem}
\Proof 
(i):
Let $0<s<\alpha\le 1$, and $\|x\|=1$, $x_1,\ldots,x_m\in \BE$. Put $t=\sum_i
\|x_i\|$. Then, by the triangle inequality and  \cite[Cor. 3.2]{Szu:hyps},
\[
\|x+Y\sum_{i=1}^m x_i\|_s \le \|1+|Yt| \|_s\le (1+c_1 t^\alpha)^{1/\alpha}
\]
for some constant $c_1=c_1(\alpha,s)$. On the other hand, since the
$\ell^1$-norm dominates the $\ell^\alpha $-norm, and by Lemma \ref{estimate},
we have
\begin{equation}\label{lowerstab}
\|x+\sum_{i=1}^m Y_i x_i\|_s\ge 
(1+c_2^p\sum_{i=1}^m \|x_i\|^\alpha)^{1/\alpha}\ge (1+C_2 t^\alpha)^{1/\alpha},
\end{equation}
which yields the assertion of the theorem, by virtue of Lemma \ref{dec-up},
where all entries of $[\bm S,{\Bbb S}]$ are equal 1.

(ii): Almost the same argument works for integrable stable chaoses.
Denote now $t=(\sum_i \|x_i\|^q)^{1/q}$, where $q\ge \alpha>1$.  By Fubini's
theorem, and the smoothness property (assuming that the norm is already
$q$-smooth), and by the hypercontractivity of Rademacher random variables, we
have
\[
\|x+Y\sum_{i=1}^m x_i\eps_i\|_s \le \|(1+|c_3Yt|^q)^{1/q} \|_s\le (1+
c_4t^\alpha)^{1/\alpha}
\]
for some constants $c_3$ and $c_4$. Since the right inequality in
\refrm{lowerstab} holds also for $\alpha>1$, we complete the proof of the
second assertion, in view of Lemma
\ref{dec-up}, with Rademacher multipliers.
\vspace{5pt}

(iii): Assume that $k^0<\infty$, i.e., $\BE$ is of Krivine's type $\le
q<\infty$.  First, let $k^0<\alpha$, and choose $q$ such that $1<q<\alpha$. By
the Figiel-Johnson renorming theorem \cite{FigJ}, there exists an equivalent
norm such that the type $\le q$-constant is equal 1.

By the hypercontractivity of the $S\alpha S$ (or $S\alpha P$) law, we may use
any $s$-norm, for $h_\alpha<s<p$, with a constant $c=c_{\alpha,q,s}(Y)$, cf.
\cite{Szu:hyps}). It is important that $h_\alpha<1$. Choose $s=q$.

We will check the following inequality  in the real case
\begin{equation}\label{auxR}
(\E|x+ X\sum_i x_i \eps_i|^q)^{1/q}
\le (\E|x+ b\sum_i x_i X_i|^q)^{1/q},
\end{equation}
where $b=b_{\alpha,q}$. Indeed, assuming that $x=1$, we obtain the following
upper bounds of the left hand side,  by virtue of the Fubini's theorem and the
hypercontractivity of Rademacher random variables ($h=h_{q,\alpha}\cdot
h_{\alpha,2}=\sqrt{(q-1)/(\alpha-1)}$,
\[
(\E(1+ |\sum_i x_i \eps_i|^{\alpha})^{q/\alpha})^{1/q}\le (\E(1+ |\sum_i x_i
\eps_i|^{\alpha})^{q/\alpha})^{1/q}\le\\ (1+h (\sum_i
|x_i|^2)^{\alpha/2})^{1/\alpha}.
\]
The right hand side is estimated from below as follows:
\[\begin{array}{c}
(\E|1+ b\sum_i x_i X_i|^q)^{1/q}\ge (1+ba_{\alpha,q}
(\Sum_i|x_i|^\alpha)^{1/\alpha}
\end{array}
\]
in view of Lemma \ref{estimate}. These estimates prove \refrm{auxR}, with
$b=h/a$.

Whence, and also by  the Fubini's theorem and hypercontractivity of $S\alpha S$
(or $S\alpha P$) law, we have
\[
\begin{array}{c}
\|x+Y\Sum_i x_i \eps_i\|_q\le 
\|(\E|x+ X\Sum_i x_i \eps_i|^q)^{1/q}\|\\
\le \|(\E|x+ b\Sum_i x_i X_i|^q)^{1/q}\|
\le \|\E|x+ bc_{q,1} \Sum_i x_i X_i|\|
\le \E\|x+ bc_{q,1} \Sum_i x_i X_i\|\\
\le \|x+ bc_{q,1} \Sum_i x_i X_i\|_q\\
\end{array}
\]
By applying Lemma \ref{dec-up}, we complete the proof of assertion (ii) in the
case $k^0<\alpha$.

Let now  $\alpha< k_0$ ($k_0\le k^0$). Choose $q\in (\alpha,k_0)$.  This case
follows immediately from assertion (ii), since
in presence of finite cotype, there is an equivalent $q$-smooth norm (cf.
\cite{Fig} or \cite[Theorem 1.f.1]{LinT}). \QED

\begin{remark}\label{SUD}\rm
Consider $S\alpha S$ (or $S\alpha P$) symmetric chaoses.
That the Krivine's classification does not fully describe the fulfillment of
the upper inequality follows from the following observation. Consider the case
$k_0\le \alpha\le k^0$.

Note that an upper decoupling inequality for $S\alpha S$ (or
$S\alpha P$) chaoses, held in Banach spaces $(\BE_1,\|\cdot \|_1)$ and
$(\BE_2,\|\cdot \|_2)$, with constants $(c_{1k})$ and $(c_{2k})$, respectively,
holds also in $\BE_1\oplus_s\BE_2$, $1\le s\alpha$, endowed with the norm
$\|\cdot\|=(\|\cdot\|_1^s+\|\cdot\|_2^s)^{1/s}$, with constants
$c_k=\max(c_{1,k},c_{2,k})$. Thus, since by assertion (ii) an upper decoupling
inequality holds in every $L^p$, $p\neq\alpha$, it will be fulfilled in every
$L^q\oplus_s L^{r}$, $q<\alpha< r$.
\end{remark}

\begin{proposition}\label{ULDstab to R}\rm
Let $\BE$ be a Banach lattice of finite cotype such that $k_0>\alpha$. Then
there exists an equivalent renorming such that all lower and upper, symmetric
and tetrahedral, decoupling constants for stable chaoses are equivalent to the
corresponding constants in the real line, i.e., $c_k(\BE)=a_{\alpha,s}^k
c_k(\R) $.
\end{proposition}

\Proof 
By a result from \cite{FigJ}, one can choose an equivalent norm of type $\ge q$
with the constant   equal to 1,  $q<\alpha $. We will use the
hypercontractivity of stable (or Pareto) (one may use any $s$-norm, for
$h_\alpha<s<\alpha $ (where $h_\alpha<1$), with a constant $a_{\alpha,s}$
\cite{Szu:hyps}). Now, denoting by $\bm Q$ and $\bm Q'$ two type of chaoses
under interest, and combining the estimates
\[
\|\Sum_k Q_k\|_s\le \|(\E|\sum_k Q_k|^s)^{1/s}\|\le \|(\E|\Sum_k (c_k(\R))^k
Q'_k|^s)^{1/s}\|
\]
and
\[
\begin{array}{c}
\|(\E|\Sum_k Q'_k|^s)^{1/s}\|\le 
\| \E|\Sum_k (a_{\alpha;s,1})^k Q'_k|\,\|\\
\le \E\|\Sum_k (a_{\alpha;s,1})^k Q'_k\| \le 
\|\Sum_k (a_{\alpha;s,1}a_{\alpha;1,s})^k Q'_k\|_s, \\
\end{array}
\]
we complete the proof.\QED
\section{CONCLUDING REMARKS}
In this section we display some further features of infinite order decoupling
and domination. Some properties or generalizations can be obtained by well
known routines, while other properties, enjoyed by homogeneous chaoses, yield
to the dead end. Yet a number of open problems arise that have no counterparts
for homogeneous chaoses. At this time, the infinite order approach to random
chaoses is still in a preliminary stage.
\subsection{Multiple stochastic integrals}
Decoupling inequalities for infinite order Gaussian or stable polynomials can
be carried over to infinite order multiple stochastic integrals, preserving all
constants, the dependence on geometry, and subjection to the presence or lack
of symmetry of underlying functions. These results follow by a routine
approximation (integrals of simple functions are random chaoses).

The real case does not require any comments, since the theory is classical. In
the vector case, one needs an appropriate construction of $k$-tuple stochastic
integrals of deterministic functions with respect to a Gaussian (or more
generally, a second order symmetric) process.  One may apply the Dunford-Bartle
approach, which reduces the integration in Banach space to that with respect to
an $L^2$-valued vector measure (cf., e.g.
\cite{DieU}).

\subsection{Non-multiplicative functions}
In \cite[Theorem 4.1]{Szu:dechom} (and before, in \cite{McCT:bandec,Pen}), a
nonmultiplicative version of the decoupling principle for homogeneous chaoses
was proved. In such a version, a term $f(i_1,\ldots,i_k)\cdot X_{1i_1}\cdots
X_{ki_k}$ was replaced by a term $F(\bm i, X_{1i_1}\cdots X_{ki_k})$.  Let us
consider a nonhomogeneous analog of such a decoupling principle (as in
\cite[4.1]{Szu:dechom}).  Let  $\L$ be an Orlicz space induced by a strongly
convex function $\varphi$ \refrm{Stass-by1}. For the sake of simplicity of
formulations, assume that $\varphi$ grows moderately.  Let $\bm F=(F_\ba)$ be a
function whose components are functions $F_\ba:\N^\ba\times \R^\ba\to \E$
satisfying conditions \cite{Szu:dechom}
\begin{equation}\label{F-conds}
\begin{array}{rl}
\mbox{\rm \hspace{-10pt}  (F1)\hspace{10pt} }&   
\mbox{\sf $F(\bm i,\cdot)=0 \quad \mu^k$-a.s. 
for all but finitely many $\bm i$;}\\
\mbox{\rm \hspace{-10pt}  (F2)\hspace{10pt} }&
\mbox{\sf $F(\bm i\,;X_{i_1},\ldots,X_{i_k})\in L^\varphi(\BE)$ for every
$\bm i\in\N^k$.}\\
\end{array}
\end{equation}

Put 
\[
\bm F(\BX^\otimes)=\sum_\ba F_\ba(\BX^{\otimes\ba}).
\]
If $\bm w=(w_\ba)$ is a Walsh sequence, write $\bm F\bm w=(F_\ba w_\ba)$ (i.e.
$[F_\ba w_\ba](\ba\bm i)=F_\ba(\ba\bm i) \bm w_\ba$).  Then the analog of
Theorem \ref{dom-by-chaos} holds, where $F_\ba$ vanish, unless $\ba=[1,k]$.
\begin{theorem}\label{F low}\rm
Let  $\L$ be an Orlicz space induced by a strongly convex function $\varphi$
\refrm{Stass-by1}, $\bm F=(F_k)$ satisfy (F1)-(F2), $\|F_k(\bx^{\otk}])\in
L^\varphi$, $k\ge 0$, and $[\bx,\BX]$ be as in Theorem \ref{dom-by-chaos}.
Then
\[
\E\varphi(\|\sum_{k\ge 0} F_k(\bx^\otk)\|)\le 
\E\varphi(\|\sum_{k\ge 0} w_k \frac{(2ck)^k}{k!}F_k(\bx^\otk)\|),
\]
where $c$ depends on the convexity of $\varphi$.
\end{theorem}
The upper decoupling inequality for functions $\bm F$ shares all deficiencies
of the corresponding decoupling inequality for homogeneous chaoses. But there
arise significant difficulties that cannot be removed by using techniques based
on hypercontractivity, since the latter method works efficiently only for
symmetric random variables. In the proof of \cite[Theorem 4.1]{Szu:dechom},
nonsymmetric random variables were used, which does not allow one to proceed as
in the proof of Theorem \ref{UDlust}. A very limited, almost trivial, real
line- version of the upper decoupling inequality can be seen as follows.  
\[
\E|\sum_k  F_k((\eps \bx)^\otk)|^2=\E|\sum_k c_k F_k(({\Bbb S} \BX)^\otk)|^2,
\]
where  $F_k(\BX^\otk)\in L^2$, and $c_k=1$ for tetrahedral functions, and
$c_k=k!$ for symmetric functions.
Any non-trivial extension (beyond Hilbert space and $L^2$-norm) would require
some intrinsic symmetry of functions $F_k$.
Therefore, at this stage it is meaningless to look at these kinds of
decoupling inequalities from the view point of integration with respect to
empirical measures (as in \cite{Szu:dechom}), even though other
types of domination might be still of interest.

\subsection{Ces\`aro averages}
There exists a variety of operators acting on the entire matrix $\BX$. For
example, one may use the operator $\bsD$, which nullifies diagonal values of
functions $f_\ba$. For the sake of consistency, denote the basic symmetrizator
by $\bsT$, $\bsT(\bm f)=\wh{\bm f}$.
Many an operator do not have meaning for a single homogeneous polynomial. We
will consider a certain multilinear analog (one of many) of Ces\`aro averages.
Let us confine ourselves to subsets $\ba\subset [1,n]$, and functions $\bm
f=(f_ba:\ba\subset[1,n])$.  We introduce the ``index average'' operator
$\bsA=(\sA_\ba)$, which unifies values of functions $ f_\ba$ along sets $\ba$
with the same cardinality.

First, we define the symmetrizator $\bsA'=(\bsA'_k)$, which transforms $\bm f$
into a function  $\bm g=(g_k:k=0,1,\dots,n)$, where
$g_k:\N^k=\N^{[0,k]}\to\BE$.

Let $|\ba|=k$. Denote by $\sfs_\ba$ the ``stretching map'' which embeds
$\N^k=\N^{[0,k]}$ into $\N^\ba$ by moving the elements of a sequence
$\bm i_k=(i_1,\ldots,i_k)=(i_1,\ldots,i_k,0,\ldots)$ into the places marked by
the consecutive ones of the sequence $\ba=(\alpha_1,\alpha_2,\ldots )$, and
filling up the remaining places by zeros. Put, for $\bm i_k=(i_1,\ldots,i_k)$,
\[
\sA'_k(\bm f)\,(\bm i_k)=\rec{{n\choose k}}
\sum_{|\alpha|=k} f_\ba(\sfs_\ba \bm i_k).
\]
Clearly,
\[
\lll\bm f \rrr=\sum_{k=0}^n{n \choose k}\ll\sA'_k(\bm f)\rr.
\]
Denote by $\sfc_{\ba}$ the ``contracting'' mapping from $\No^\ba$ onto $\No^k$,
which just cancels all elements marked by zeros of the sequence $\ba$. Now, we
define the ``inverse'' mapping $\bsA''=(\sA''_\ba)$ transforming functions
$\bm g=(g_k)$  into functions $\bm f=(f_\ba)$, according  to the formula
\[
\sA''_\ba(\bm g)\,(\ba\bm i)=g_k(\sfc_\ba(\ba\bm i)),    \qquad|\ba|=k.
\]
Define $\bsA=\bsA''\bsA'$. 
The operators $\bsD$, $\bsT$, and  $\bsA$ are idempotent and  commute
with each other. 

Let $\BE_1,\BE_2$ be additive abelian groups. Denote by $x_1x_2$
a bi-additive mapping from $\BE_1\times \BE_2$ into $\BE$. Use the same
notation $\bm f_1\bm f_2$ for functions taking values in $\BE_1$ and $\BE_2$,
respectively. If $\bsU$ and $\bsV$ are compositions of selected 
symmetrizators $\bsD$, $\bsA$, $\bsT$, then the following 
symmetrization formulas hold:
\begin{equation}\label{symUV}
\lll \bsU(\bm f_1)\bsV(\bm f_2) \rrr =
\lll \bm f_1\bm{\bsU\bsV}(\bm f_2) \rrr =
\lll \bsU\bsV(\bm f_1)\bm f_2 \rrr =
\lll \bsV(\bm f_1)\bsU(\bm f_2) \rrr. 
\end{equation}
Note that $\bx^\otimes$ is $\bsA$-symmetric.
By $\ca A_k=\ca A_k(\BX)$ denote the $\sigma$-field generated by the family of
random variables
\[
\set{h(\BX^\ba:|\ba|=k)\,:\,h=\hat{h}, \,h:(\R^k)^{n\choose k} \to \R}.
\]

Notice that the symmetry assumption is applied to $h$ as to a function of ${n
\choose k }$ vector variables, and that $\ca A_k$ are ascending
$\sigma$-fields. The symmetrizator $\bsA$  can be expressed as  a conditional
expectation. Note the following equalities:

\begin{equation}\label{Acontr}
\lll \bm f \BX^{\bsA(\otimes)} \rrr =
\lll \bm f \E[\BX^{\otimes}|\ca A({\BX})] \rrr =
\E[\lll \bm f \BX^{\otimes}\rrr|\ca A(\BX), 
\end{equation}
\begin{equation}\label{Ave}
\E\left[\, \bsD(\bx_1+\ldots+\bx_k)^{\otimes k}\,
|\,\ca A(\BX) \,\right]
=\bsD \sA_\ba\l(\ba\BX\r)^{\otimes\ba};
\end{equation}
or equivalently,
\begin{equation}\label{Avesum}
\E\left[\,\bsD(\frac{\bx_1+\ldots+\bx_k}{k})^{\otimes k}\,|\,\ca A(\BX)
\,\right] =\bsD A_\ba\left(\frac{\ba\bx}{|\ba|}\right)^{\otimes\ba}.
\end{equation}
Now, Theorem \ref{dom-by-chaos} holds for $\bsA$-convex functions. That is, 
\begin{equation}\label{A dom-by-chaos}
\E\varphi(\|\San w_\ba\ll  f_\ba \BX^{\otimes\ba} \rr \|)\le 
\E\varphi(\|\San w_\ba h_\ba\ll  f_\ba \bx^{\otimes\ba} \rr \|), 
\end{equation}
for every $\bsD-,\bsA, \&\bsT$-symmetric function $\bm f$, where, for
$|\ba|=k$, $h_\ba=h_k=(2ck)^k/k!$, and $c=c_\varphi$.

However, the $\bsA$-symmetry is too strong for the upper  inequality of
arbitrary order to be fulfilled. For integrable symmetric random variables, by
examining just polynomials of the first degree, we would obtain the inequality
\[
\E\varphi (\|x+ \sum_{i\le K} X_i x_i \|) \le 
\E\varphi (\|x+ c_1 \sum_{i\le K} \frac{\sum_{j=1}^n X_{ji}}{n} x_i\|) 
\]
which is impossible, as can be seen by applying the strong law of large numbers
and Fatou's lemma. Yet, the above observations open a new, even in the real
case, direction in a study of $\bsA$-symmetric chaoses. Clearly, any domination
``constant'' is expected to depend on $n$, which makes the concept of infinite
order  much more difficult.

{\bf Problem}. Describe the closure in $L^2$ and limit distributions of real
Gaussian (or Rademacher) $\bsA$-symmetric decoupled chaoses
\[
\set{\sum_{\ba\subset[1,n]} \ll f_\ba \BX^{\bsA(\ba)}\rr
:n\in\N}.   
\]
Note that the metric ($L^2$-) problem is easy for $\bsT$-symmetric or
tetrahedral functions (cf. the first subsection of this section). For
$\bsT$-symmetric functions, a related limit theorem for coupled Gaussian random
chaoses, obtained in \cite{DynM}, brought up infinite order Wiener integrals.

\renewcommand{\baselinestretch}{0.8}
\sf
\vfill

~\hfill 
\begin{tabular}{|l|}\hline
\small Jerzy Szulga\\
\small Department of Mathematics, 
       Mathematics Annex 120\\
\small Auburn University, 
       Auburn, AL 36849-3501\\
\begin{tabular}{rl}
\small phone &\small (205) 844-3649\\
\small email &szulga@auducvax.bitnet\\
             &szulga@ducvax.auburn.edu\\
\end{tabular}\\ \hline
\end{tabular}
\end{document}